\documentclass[a4paper, 11pt]{article}
\usepackage{amsmath}
\usepackage{amsfonts}
\usepackage{amssymb}
\usepackage[english]{babel}
\usepackage{graphicx}
\usepackage{amsthm,cases}
\usepackage[numbers,sort&compress]{natbib}
\usepackage[title]{appendix}

\usepackage{longtable}
\usepackage{tabularx}
\usepackage{diagbox,ulem}
\allowdisplaybreaks

\newcommand{\B}{{\cal B}}

\newcommand{\C}{{\cal C}}
\newcommand{\F}{{\cal F}}
\newcommand{\A}{{\cal A}}

\newcommand{\G}{{\cal G}}

\newtheorem{theorem}{Theorem}[section]

\newtheorem{lemma}[theorem]{Lemma}

\newtheorem{construction}[theorem]{Construction}

\def\whitebox{{\hbox{\hskip 1pt
 \vrule height 6pt depth 1.5pt
 \lower 1.5pt\vbox to 7.5pt{\hrule width
    3.2pt\vfill\hrule width 3.2pt}%
 \vrule height 6pt depth 1.5pt
 \hskip 1pt } }}
\def\qed{\ifhmode\allowbreak\else\nobreak\fi\hfill\quad\nobreak
     \whitebox\medbreak}

\newcommand{\ignore}[1]{}
\allowdisplaybreaks

\begin{document}

\title{ The spectrum for  large sets of $(3,\lambda)$-GDDs of type $g^u$}

\author{\small  X. Niu $^1$, \ H. Cao $^1$ \thanks{Research
supported by the National Natural Science Foundation of China
under Grant 11571179, and the Priority Academic
Program Development of Jiangsu Higher Education Institutions. E-mail: caohaitao@njnu.edu.cn
} ,\  and R.  Javadi $^{2,~3}$ \thanks{The research is partially carried out in the IPM-Isfahan Branch and
in part supported, respectively, by grant No. 96050059 from IPM. E-mail: rjavadi@cc.iut.ac.ir
}\\
\small $^1$ Institute of Mathematics,  Nanjing Normal
University, Nanjing 210023, China\\
\small $^2$ Department of Mathematical Sciences,  Isfahan University of Technology, Isfahan 84156-83111, Iran\\
\small $^3$ School of Mathematics, Institute for Research in Fundamental Sciences, Tehran 19395-5746, Iran\\}

\date{}
\maketitle
\begin{abstract} In this paper, we completely solve the existence of large sets of $(3,\lambda)$-GDDs of type $g^u$ and the existence of a simple $(3,\lambda)$-GDD of type $g^u$.

\bigskip



\end{abstract}

\section{Introduction}

A {\it group divisible design}
GDD$(\lambda, t,K,v;g_1^{u_1}, \ldots, g_s^{u_s})$,
$v=\sum_{i=1}^{s}g_iu_i, K\subset \mathbb{N}$,
 is a triple $(X,\G,\B)$ such that

(1) $X$ is a set of $v$ elements (called {\it points}).

(2) $\G$ is a partition of $X$ into $u_i$ sets of $g_i$ points (called {\it groups}), $i = 1,2, \ldots, s$.

(3) $\B$ is a multiset of subsets of $X$ (called {\it blocks}), such that $|B| \in K$, $|B \cap G| \le 1 $ for all
$B \in  \B$ and $G \in \G$ and such that any $t$-subset $T$ of $X$ with $|T \cap G| \le 1 $ for all $G \in \G$,
is contained in exactly $\lambda$ blocks.

For the sake of brevity, we write $k$ rather than $\{k\}$ when $K=\{k\}$ is a singleton,
and use  $(k,\lambda)$-GDD$(g_1^{u_1}, \ldots, g_s^{u_s})$
instead of  GDD$(\lambda,2,k,v;g_1^{u_1}, \ldots, g_s^{u_s})$.
A GDD$(\lambda,t,K,v;1^{v})$ is called a {\it t-wise balanced design},
denoted by $S(\lambda;t,K,v)$.
If $\lambda=1$, then $\lambda$ is usually omitted. An $S(t,k,v)$ is called a {\it Steiner system}. An  S$(\lambda;2,  K,  v)$  is called a {\it  pairwise balanced design} or PBD.

A  group divisible design is called {\it simple} if no two blocks are identical.
In this paper, we focus on group divisible designs with block size 3.
The necessary condition for the existence of a simple $(3,\lambda)$-GDD$(g^u)$ is as follows.

\begin{theorem}\label{Ne-SGDD}
If there is a simple $(3,\lambda)$-GDD$(g^u)$, then $u\ge 3$, $1 \leq \lambda \leq g(u-2)$,
 $\lambda g(u-1)\equiv 0 \pmod 2$, and $\lambda g^2u(u-1)\equiv 0 \pmod 6$.
\end{theorem}

Two simple $(3, \lambda)$-GDDs $(X,\G,\A)$ and $(X,\G,\B)$ are called {\it disjoint}
if ${\A}\cap {\B}=\emptyset$. A set of more than two simple $(3, \lambda)$-GDDs is called {\it disjoint}
if each pair of them is disjoint. It is obvious that the maximum number of disjoint $(3, \lambda)$-GDDs of type $g^u$ is no more than $\frac{(u-2)g}{\lambda}$.
The collection  of $\frac{(u-2)g}{\lambda}$ disjoint $(3, \lambda)$-GDDs of type $g^u$ is called a {\it large set}. We denote the large set by  $(3,\lambda)$-LGDD$(g^u)$. The necessary condition for the existence of a $(3,\lambda)$-LGDD$(g^u)$ is as follows.

\begin{theorem}\label{Ne-LGDD}
If there is
a $(3,\lambda)$-LGDD$(g^u)$, then  $u\ge 3$,  $1 \leq \lambda \leq g(u-2)$,
$\lambda g(u-1)\equiv 0 \pmod 2$, $\lambda g^2u(u-1)\equiv 0 \pmod 6$ and $g(u-2)\equiv 0 \pmod {\lambda}$.
\end{theorem}

A lot of work has been done on the existence of a simple  $(3,\lambda)$-GDD$(g^u)$ and a $(3,\lambda)$-LGDD$(g^u)$. We summarize the known results in the following two theorems. For more results on large sets with resolvable property, see \cite{CZ2013,JL2008,L2002,T1994,ZCZ2017,ZCZ2017-1}.

\begin{theorem}\label{SGDD}{\rm (\cite{D1983,S1974,T1977,VB1974-1,VB1974-2,W1974})}
There exists a simple $(3,\lambda)$-GDD$(1^u)$
if and only if $1\le \lambda\le u-2$, $\lambda(u-1)\equiv 0 \pmod 2$ and $\lambda u(u-1)\equiv 0 \pmod 6$.
\end{theorem}

\begin{theorem}\label{LGDD}
1. There exists a $(3,1)$-LGDD$(g^u)$ if and only if $g(u-1)\equiv 0 \pmod 2$, $g^2u(u-1)\equiv 0 \pmod 6$, $u\ge 3$ and $(g,u)\neq(1,7)$ {\rm \cite{CLS,CS1,CS,L1997,Lu1,Lu2,SS1989,SV1988,Teir,T1993}}.

\noindent 2. There exists a $(3,2)$-LGDD$(1^u)$ if and only if  $u\equiv 0, 4 \pmod 6$ and $u\ge 4$ {\rm \cite{S1974,T1975}}.

\noindent 3. There exists a $(3,3)$-LGDD$(1^u)$ if and only if   $u\equiv 5 \pmod 6$ and $u\ge 5$ {\rm \cite{T1984}}.

\noindent 4. There exists a $(3,6)$-LGDD$(1^u)$ if and only if   $u\equiv 2 \pmod 6$ and $u\ge 8$ {\rm \cite{T1975}}.
\end{theorem}

By definitions it is easy to see that a $(3,\lambda)$-LGDD$(g^u)$ can be used to obtain a simple $(3, t\lambda)$-GDD$(g^u)$ for any $1\le t\le \frac{g(u-2)}{\lambda}$. So we have the following lemma.

\begin{lemma}\label{LS-SGDD}
If there exists a $(3,\lambda)$-LGDD$(g^u)$, then there exists a simple $(3,t\lambda)$-GDD$(g^u)$ for any $1\le t\le \frac{g(u-2)}{\lambda}$.
\end{lemma}

In this paper, we shall obtain necessary and sufficient conditions for the existence of a $(3,\lambda)$-LGDD$(g^u)$ and the existence of a simple $(3,\lambda)$-GDD$(g^u)$ for any $g>1$.
Generalizing the known results in Theorems~\ref{SGDD} and \ref{LGDD}, we
will prove the following two theorems.

\begin{theorem}\label{Main-L}
There exists a $(3,\lambda)$-LGDD$(g^u)$ if and only if
$u\ge 3$, $1 \leq \lambda \leq g(u-2)$, $\lambda g(u-1)\equiv 0 \pmod 2$,
$\lambda g^2u(u-1)\equiv 0 \pmod 6$, $gu(u-2)\equiv 0 \pmod {\lambda}$ and $(\lambda,g,u)\neq (1,1,7)$.
\end{theorem}

\begin{theorem}\label{Main-S}

There exists a simple $(3,\lambda)$-GDD$(g^u)$ if and only if
$u\ge 3$, $1 \leq \lambda \leq g(u-2)$,
 $\lambda g(u-1)\equiv 0 \pmod 2$ and $\lambda g^2u(u-1)\equiv 0 \pmod 6$.
\end{theorem}

This paper is organized as follows. In the next section, we will introduce some necessary definitions and notations, and generalize L. Teirlinck's method called large sets with holes which will play an important role in our constructions.
In Section 3, we shall give some constructions for some good large sets defined in Section 2 by using LR designs and generalized frames.
In the last two sections, we shall prove Theorems~\ref{Main-L} and \ref{Main-S}.

\section{Large sets with holes  }

In this section,
we shall introduce a method called large sets with holes which was posed by L. Teirlinck in \cite{T1993},
and define some new large sets.

An LS$(\lambda_1,\lambda_2;2,(3,K),v)$ is a
collection  $(X, \B_r )_{r\in R}$ of S$(\lambda_2;2,K,v)$ such that
$(X, \cup_{r\in R} \B_r)$
 is an S$(3,K,v)$ and such that,
for each $B\in \cup_{r\in R} \B_r$,
$B$ appears  exactly $\lambda_1(|B|-2)$ times in the multiset $\{B: B\in\B_r, r\in R\}$. (Note that $\cup_{r\in R} \B_r$ denotes the ordinary union of the $\B_r$ and not the multiset-union.)
 L. Teirlinck mainly considered the special case LS$(\lambda, 1;2,(3,K),v)$ in \cite{T1993}.
In this paper, we will focus on another special case LS$(1,\lambda;2,(3,K),v)$. We usually  write LS$(2,(3,K),v)$ instead of LS$(1,1;2,(3,K),v)$.
An LS$(1,\lambda; 2,(3,\{3\}),v)$ is essentially equivalent to a $(3,\lambda)$-LGDD$(1^v)$.

We first define a new large set with special properties which will play an important role in our recursive constructions.
An  LS*$(1,\lambda;2,(3,K),v)$ is an LS$(1,\lambda;2,(3,K),v)$  $(X, \B_r )_{r\in R}$ such that
 for any $r\in R$, $B\in \B_r$ and $|B|\ge 4$,
$B$ appears exactly $\lambda$ times in the $\B_r$.

An S$(\lambda; 2, (K_{0}, K_1, K_2), v)$ is a quadruple $(S, \infty_{1}, \infty_{2}, {\cal B})$,
where  $|S|=v-2$, $S\cap\{\infty_{1}, \infty_{2}\}= \emptyset$, $\infty_1\not=\infty_2$,
such that $(S\cup \{\infty_{1}, \infty_{2}\}, {\cal B})$ is an S$(\lambda; 2, K_0\cup K_1\cup K_2, v)$
satisfying  $|B|\in K_{i}$  for all  $B\in {\cal B}$  with  $|B\cap \{\infty_{1}, \infty_{2}\}| =  i$, $i=0,1,2$.
An LS$(1, \lambda; 2, (3, K_{0}, K_1, K_2), v)$
will be a collection $(S, \infty_{1}, \infty_{2}, {\cal B}_{r})_{r\in R}$
of S$(2, (K_{0}, $ $K_1, $ $K_2), v)$ such that
$(S\cup \{\infty_{1}, \infty_{2}\}, {\cal B}_{r})_{r\in R}$
is an LS$(1,\lambda;2,(3, K_{0}\cup K_1\cup K_2),v)$.
In fact, an LS$(1,\lambda; 2, $ $(3, K_{0}, K_1, K_2), v)$ is the same thing as
 an   LS$(1,\lambda; 2, (3, K_{0}\cup K_1\cup K_2), v)$ with two distinguished points.
The proofs of the following constructions for LS, LS* and LGDD are similar to the proofs of Constructions 2.2$(b)$, 4.1, and 5.1 in \cite{T1993}.

\begin{construction}\label{w=2}
\noindent (1) If there exist an LS$(1,\lambda; 2, (3, K_0, \{3\}, K_2), v)$,
 an LS$( 2,(3,K_0^{\prime}), k)$ for all $k\in K_0$, and
 an LS$(2,(3,K_0^{\prime},\{3\}, K_2^{\prime} ), k)$ for all $k\in K_2$ ,
 then there exists an LS$(1,\lambda; 2,(3,K_0^{\prime},\{3\},K_2^{\prime}), v)$.

\noindent (2) If there exist an LS$(1,\lambda; 2, (3, K_0, \{3\}, K_2), v)$ and
an $(3,1)$-LGDD$(2^k)$ for all $k\in K_0$,
then there exists an
LS$(1,\lambda; 2,(3,\{3, 4\},\{3\}, K_{2}^{\prime}),2v-2)$,
where $K_{2}^{\prime} = \{2k-2 : k\in K_2\}$.

\noindent (3) If there exist an LS$(1,\lambda;2,(3,K),v)$
and  a $(3,1)$-LGDD$(g^k)$ for all $k\in K$,
then there exists a $(3,\lambda)$-LGDD$(g^v)$.
\end{construction}

\begin{construction}\label{LS-LGDD}
\noindent (1) If there exist an LS*$(1,\lambda; 2, (3, K_0, \{3\}, K_2), v)$,
an LS*$(1,\lambda; 2,(3,K_0^{\prime}),$ $ k)$ for all $k\in K_0$, and
 an LS*$(1,\lambda;2,(3,K_0^{\prime},\{3\}, K_2^{\prime} ), k)$ for all $k\in K_2$,
 then there exists an LS*$(1,\lambda; 2,(3,K_0^{\prime},$ $\{3\},K_2^{\prime}), v)$.

\noindent (2) If there exist an LS*$(1,\lambda; 2, (3, K_0, \{3\}, K_2), v)$,
an LS$(1,\lambda; 2,(3,K_0^{\prime}),$ $ k)$ for all $k\in K_0$, and
 an LS$(1,\lambda;2,(3,K_0^{\prime},\{3\}, K_2^{\prime} ), k)$ for all $k\in K_2$,
 then there exists an LS$(1,\lambda; 2,(3,K_0^{\prime},$ $\{3\},K_2^{\prime}), v)$.

 \noindent (3) If there exist an LS$( 2, (3, K_0, \{3\}, K_2), v)$ and
 an $(3,1)$-LGDD$(2^k)$ for all $k\in K_0$,
then there exists an LS*$(1,2; 2,(3,$ $\{3, 4\},\{3\}, K_{2}^{\prime}),2v-2)$,
where $K_{2}^{\prime} = \{2k-2 : k\in K_2\}$.

\noindent (4) If there exist an LS*$(1,\lambda;2,(3,K),v)$
and  a $(3,\lambda)$-LGDD$(g^k)$ for all $k\in K$ and $k\ge 4$,
then there exists a $(3,\lambda)$-LGDD$(g^v)$.
\end{construction}

Now we continue to introduce the concept of good large sets defined by L. Teirlinck in \cite{T1993}.
A{\it digraph} is an ordered pair $D = (V(D),A(D))$,
where $V(D)$ is a set whose elements are called {\it vertices}
and $A(D)$ is a set of ordered pairs of vertices, called {\it directed edges}.
A directed edge $(x, y)\in A(D)$, $x,y\in V(D)$, is considered to be directed from $x$ to $y$.
$y$ is called the {\it head} and $x$ is called the {\it tail} of the directed edge $(x, y)$.
For a vertex $x\in V(D)$, the {\it indegree} $d^{-}(x)$ in $D$ is the number of directed edges with head $x$,
and the {\it outdegree} $d^{+}(x)$ is the number of directed edges with tail $x$.
If for every vertex $x\in V(D)$, $d^+(x) = d^-(x)$,
the graph $D$ is called an {\it Eulerian digraph}.
It is easy to check that the union of any two Eulerian digraphs is also an Eulerian digraph.

A {\it good} S$(\lambda; 2,(K_{0}, \{3\}, K_2), v)$ (or GS$(\lambda; 2,(K_{0}, \{3\}, K_2), v)$)
is a 5-tuple $(S, \infty_{1}, \infty_{2}, {\cal B}, D)$,
where  $(S, \infty_{1}, \infty_{2}, {\cal B})$ is an S$(\lambda; 2, (K_{0}, \{3\}, K_2), v)$,
and where $D$ is an Eulerian digraph on $S$ whose underlying undirected graph has edge-set
$\{\{x, y\} : x, y  \in S, \{\infty_{i}, x, y\} \in {\cal B}, i \in \{1,2\}\}$.
Let ${\cal A}_i=\{B : B \in{\cal B}, B\cap \{\infty_{1}, \infty_{2}\}=\{\infty_{i}\} \}$, $i= 1,2$.
For any  $x\in S$, let $t_x=|\{B : \{\infty_{1}, \infty_{2}, x\} \subset B, B \in{\cal B} \}|$,
 $d_i^+(x)=|\{(x,y): y\in S, (x,y)\in A(D), \{\infty_{i}, x, y\} \in {\cal A}_i \}|$ and
 $d_i^-(x)=|\{(y,x): y\in S, (y,x)\in A(D), \{\infty_{i}, x, y\} \in {\cal A}_i \}|$ for $i= 1,2$.
Since $(S, \infty_{1}, \infty_{2}, {\cal B})$ is an S$(\lambda; 2, (K_{0}, \{3\}, K_2), v)$,
then we have $d_i^+(x) +d_i^-(x)+t_x=\lambda$ for $i= 1,2$.
Since $D$ is an Eulerian digraph on $S$,
then we have $d_1^+(x) +d_2^+(x)=d^+(x)=d^-(x)=d_1^-(x) +d_2^-(x)$.
Thus, $d^+(x)=d^-(x)=\lambda-t_x$.

A {\it good} LS$( 1,\lambda; 2, (3, K_0, \{3\}, K_2), v)$
(or GLS$(1,\lambda; 2, (3, K_0, \{3\}, K_2), v)$) will be a collection
$(S, \infty_{1}, \infty_{2}, {\cal B}_r, D_r)_{r\in R}$ of GS$(\lambda; 2,(K_{0}, \{3\}, K_2), v)$,
where $(S, \infty_{1}, \infty_{2}, {\cal B}_r)_{r\in R}$ is an LS$(1,\lambda;2,$ $(3,K_{0}, \{3\}, K_2), v)$,
such that each ordered pair $(x, y)$ of distinct elements of $S$,
not contained in some block $B$ with $\{\infty_{1},\infty_{2}\}\subset B$, appears in exactly one $D_r$.

A GLS*$( 1,\lambda; 2, (3, K_0, \{3\}, K_2), v)$ is an LS$( 1,\lambda; 2, (3, K_0, \{3\}, K_2), v)$ satisfying all the requirements of a GLS$(1,\lambda; 2, (3, K_0, \{3\}, K_2), v)$
and an LS*$( 1,\lambda; 2, (3, K_0, \{3\}, K_2), v)$.
If $\lambda=1$,
we often write  GLS$(2, (3, K_0, \{3\}, K_2), v)$ instead of GLS$(1,1;2, (3, K_0, \{3\}, K_2), v)$
and write GLS$(2,3, v)$ instead of GLS$(1, 1; 2, (3, \{3\}, \{3\}, $ $\{3\}), v)$.
A GLS$(2,3, v)$ is essentially equivalent to a GLS$(v)$, as defined by Lu in \cite{Lu1}.
We need the following result on GLS$(2,(3, K_0, \{3\}, K_2), v)$ for our constructions.

\begin{theorem}\label{GLS} {\rm  (\cite{T1993})}
1. There exists a GLS$(2,(3,\{3\},\{3\},\{5\}), v)$ for $v\equiv 5\pmod 6$.

\noindent 2. There exists a GLS$(2,(3, \{3,4\}, \{3\}, \{8,14,26,50\}), v)$ for all $v\equiv 2 \pmod 6$.
\end{theorem}

\begin{lemma}\label{v=5}
 There exists a GLS$(1,3;2,(3, \{3\}, \{3\}, \{3\}), 5)$.
\end{lemma}
\noindent {\it Proof:}
Let $S = \{0,1,2\}$ and

{\small
\begin{tabular}{lllllllllll}
 $\B = $
&$\{( \infty_{1}, \infty_{2},0)$,
&$  (\infty_{1}, \infty_{2},1)$,
&$  (\infty_{1}, \infty_{2},2)$,
&$  (\infty_{1}, 0, 1)$,
&$  (\infty_{1}, 0, 2)$,\\
&$  (\infty_{1}, 1, 2)$,
&$  (\infty_{2}, 1, 0)$,
&$  (\infty_{2}, 2, 0)$,
&$  (\infty_{2}, 2, 1)$,
&$  (0, 1, 2)$\},\\
\end{tabular}

$A(D)  = \ \  \{(0,1),(0,2),(1,2),(1,0),(2,0),(2,1)\}$.
}

It is easy to check that $(S,\infty_{1}, \infty_{2}, \B, D)$
is a GLS$(1,3;2,(3, \{3\}, \{3\}, \{3\}), 5)$.
\qed

\begin{lemma}\label{v=6}
 There exists a GLS$(1,2;2,(3, \{3\}, \{3\}, \{3\}), 6)$.
\end{lemma}
\noindent {\it Proof:}
Let $S = \{0,1,2,3\}$, $R = \{0,1\}$, and

{\small
\begin{tabular}{lllllllllll}
 $\B_0 = $
&$\{( \infty_{1}, \infty_{2},0)$,
&$  (\infty_{1}, \infty_{2},1)$,
&$  (\infty_{1}, 0, 2)$,
&$  (\infty_{1}, 1, 3)$,
&$  (\infty_{1}, 2, 3)$,\\
&$  (\infty_{2}, 3, 0)$,
&$  (\infty_{2}, 2, 1)$,
&$  (\infty_{2}, 3, 2)$,
&$  (0, 1, 2)$,
&$  (0, 1, 3)$\},\\
\end{tabular}

\begin{tabular}{lllllllllll}
 $\B_1 = $
&$\{( \infty_{1}, \infty_{2}, 2)$,
&$  (\infty_{1}, \infty_{2}, 3)$,
&$  (\infty_{1}, 0, 1)$,
&$  (\infty_{1}, 0, 3)$,
&$  (\infty_{1}, 1, 2)$,\\
&$  (\infty_{2}, 1, 0)$,
&$  (\infty_{2}, 2, 0)$,
&$  (\infty_{2}, 3, 1)$,
&$  (0, 2, 3)$,
&$  (1, 2, 3)$\},\\
\end{tabular}

 $A(D_0) = \ \  \{(0,2),(1,3),(2,3),(3,0),(2,1),(3,2)\}$,

 $A(D_1)  = \ \  \{(0,1),(0,3),(1,2),(1,0),(2,0),(3,1)\}$.}

It is easy to check that $(S, \infty_{1}, \infty_{2}, $ $\B_r, D_r)_{r\in R}$
is a GLS$(1,2;2,(3, \{3\}, \{3\},$ $ \{3\}), 6)$.
\qed

\begin{lemma}\label{v=10}
 There exists a GLS$(1,2;2,(3, \{3\}, \{3\}, \{3\}), 10)$.
\end{lemma}
\noindent {\it Proof:}
Let $S = \{0,1,2,\ldots, 7\}$, $R=\{0,1,2,3\}$, and

{\small
\noindent\begin{tabular}{lllllllllll}
 $\B_0 = $
&$\{( \infty_{1}, \infty_{2}, 0)$,
&$  (\infty_{1}, \infty_{2}, 1)$,
&$  (\infty_{1}, 0, 4)$,
&$  (\infty_{1}, 1, 5)$,
&$  (\infty_{1}, 2, 4)$,
&$  (\infty_{1}, 2, 6)$,\\
&$  (\infty_{1}, 3, 5)$,
&$  (\infty_{1}, 3, 7)$,
&$  (\infty_{1}, 6, 7)$,
&$  (\infty_{2}, 4, 1)$,
&$  (\infty_{2}, 4, 3)$,
&$  (\infty_{2}, 5, 0)$,\\
&$  (\infty_{2}, 5, 2)$,
&$  (\infty_{2}, 6, 3)$,
&$  (\infty_{2}, 7, 2)$,
&$  (\infty_{2}, 7, 6)$,
&$  (0, 1, 2)$,
&$  (0, 1, 3)$,\\
&$  (0, 2, 6)$,
&$  (0, 3, 7)$,
&$  (0, 4, 6)$,
&$  (0, 5, 7)$,
&$  (1, 2, 7)$,
&$  (1, 3, 6)$,\\
&$  (1, 4, 7)$,
&$  (1, 5, 6)$,
&$  (2, 3, 4)$,
&$  (2, 3, 5)$,
&$  (4, 5, 6)$,
&$  (4, 5, 7)$\},\\
\end{tabular}

\noindent\begin{tabular}{lllllllllll}
 $\B_1 = $
&$\{( \infty_{1}, \infty_{2}, 2)$,
&$  (\infty_{1}, \infty_{2}, 3)$,
&$  (\infty_{1}, 0, 2)$,
&$  (\infty_{1}, 0, 6)$,
&$  (\infty_{1}, 1, 3)$,
&$  (\infty_{1}, 1, 7)$,\\
&$  (\infty_{1}, 4, 5)$,
&$  (\infty_{1}, 4, 6)$,
&$  (\infty_{1}, 5, 7)$,
&$  (\infty_{2}, 2, 1)$,
&$  (\infty_{2}, 3, 0)$,
&$  (\infty_{2}, 5, 4)$,\\
&$  (\infty_{2}, 6, 1)$,
&$  (\infty_{2}, 6, 5)$,
&$  (\infty_{2}, 7, 0)$,
&$  (\infty_{2}, 7, 4)$,
&$  (0, 1, 4)$,
&$  (0, 1, 5)$,\\
&$  (0, 2, 4)$,
&$  (0, 3, 5)$,
&$  (0, 6, 7)$,
&$  (1, 2, 5)$,
&$  (1, 3, 4)$,
&$  (1, 6, 7)$,\\
&$  (2, 3, 6)$,
&$  (2, 3, 7)$,
&$  (2, 4, 6)$,
&$  (2, 5, 7)$,
&$  (3, 4, 7)$,
&$  (3, 5, 6)$\},\\
\end{tabular}

\noindent\begin{tabular}{lllllllllll}
 $\B_2 = $
&$\{( \infty_{1}, \infty_{2}, 4)$,
&$  (\infty_{1}, \infty_{2}, 5)$,
&$  (\infty_{1}, 0, 3)$,
&$  (\infty_{1}, 0, 7)$,
&$  (\infty_{1}, 1, 2)$,
&$  (\infty_{1}, 1, 6)$,\\
&$  (\infty_{1}, 2, 3)$,
&$  (\infty_{1}, 4, 7)$,
&$  (\infty_{1}, 5, 6)$,
&$  (\infty_{2}, 2, 0)$,
&$  (\infty_{2}, 3, 1)$,
&$  (\infty_{2}, 3, 2)$,\\
&$  (\infty_{2}, 6, 0)$,
&$  (\infty_{2}, 6, 4)$,
&$  (\infty_{2}, 7, 1)$,
&$  (\infty_{2}, 7, 5)$,
&$  (0, 1, 6)$,
&$  (0, 1, 7)$,\\
&$  (0, 2, 5)$,
&$  (0, 3, 4)$,
&$  (0, 4, 5)$,
&$  (1, 2, 4)$,
&$  (1, 3, 5)$,
&$  (1, 4, 5)$,\\
&$  (2, 4, 7)$,
&$  (2, 5, 6)$,
&$  (2, 6, 7)$,
&$  (3, 4, 6)$,
&$  (3, 5, 7)$,
&$  (3, 6, 7)$\},\\
\end{tabular}

\noindent\begin{tabular}{lllllllllll}
 $\B_3 = $
&$\{( \infty_{1}, \infty_{2}, 6)$,
&$  (\infty_{1}, \infty_{2}, 7)$,
&$  (\infty_{1}, 0, 1)$,
&$  (\infty_{1}, 0, 5)$,
&$  (\infty_{1}, 1, 4)$,
&$  (\infty_{1}, 2, 5)$,\\
&$  (\infty_{1}, 2, 7)$,
&$  (\infty_{1}, 3, 4)$,
&$  (\infty_{1}, 3, 6)$,
&$  (\infty_{2}, 1, 0)$,
&$  (\infty_{2}, 4, 0)$,
&$  (\infty_{2}, 4, 2)$,\\
&$  (\infty_{2}, 5, 1)$,
&$  (\infty_{2}, 5, 3)$,
&$  (\infty_{2}, 6, 2)$,
&$  (\infty_{2}, 7, 3)$,
&$  (0, 2, 3)$,
&$  (0, 2, 7)$,\\
&$  (0, 3, 6)$,
&$  (0, 4, 7)$,
&$  (0, 5, 6)$,
&$  (1, 2, 3)$,
&$  (1, 2, 6)$,
&$  (1, 3, 7)$,\\
&$  (1, 4, 6)$,
&$  (1, 5, 7)$,
&$  (2, 4, 5)$,
&$  (3, 4, 5)$,
&$  (4, 6, 7)$,
&$  (5, 6, 7)$\},\\
\end{tabular}

\noindent $A(D_0)  =   \{(0,4), (1,5), (2,4), (2,6), (3,5), (3,7), (6,7)$,
                              $(4,1), (4,3), (5,0), (5,2), (6,3), (7,2), (7,6)\}$,

\noindent $A(D_1)  =   \{(0,2), (0,6), (1,3), (1,7), (4,5), (4,6), (5,7)$,
                              $(2,1), (3,0), (5,4), (6,1), (6,5), (7,0), (7,4)\}$,

\noindent $A(D_2)  =  \{(0,3), (0,7), (1,2), (1,6), (2,3), (4,7), (5,6)$,
                               $(2,0),(3,1),(3,2),(6,0),(6,4),(7,1),(7,5)\}$,

\noindent $A(D_3)  =   \{(0,1), (0,5), (1,4), (2,5), (2,7), (3,4), (3,6)$,
                             $(1,0),(4,0),(4,2),(5,1),(5,3),(6,2),(7,3)\}$.
}

It is easy to check that $(S,\infty_{1}, \infty_{2}, \B_r, D_r)_{r\in R}$
is a GLS$(1,2;2,(3, \{3\}, \{3\}, \{3\}), 10)$.
\qed

\begin{lemma}\label{v=11}
 There exists a GLS*$(1,3; 2,(3,\{3\},\{3\},\{5\}),11)$.
\end{lemma}

\noindent {\it Proof:}
Let $S=\mathbb{Z}_9$, $R=\{0, 1,2\}$, and

{\small
\begin{tabular}{lllllllllll}
 $\A_0 = \{$
&$  (\infty_{1}, 3i +1, 3i+5)$,
&$  (\infty_{1}, 3i +4, 3i+2)$,
&$  (\infty_{1}, 3i +7, 3i+8)$,\\
&$  (\infty_{2}, 3i+8, 3i+1)$,
&$  (\infty_{2}, 3i+5, 3i+4)$,
&$  (\infty_{2}, 3i+2, 3i+7)$,\\
&$  (3i, 3i+1, 3i+2)$,
&$  (3i, 3i+4, 3i+8)$,
&$  (3i, 3i+7, 3i+5)$,\\
&$  (3i, 3i+4, 3i+7)$,
&$  (3i, 3i+7, 3i+1)$,
&$  (3i, 3i+1, 3i+4)$,\\
&$  (3i, 3i+5, 3i+8)$,
&$  (3i, 3i+2, 3i+5)$,
&$  (3i, 3i+8, 3i+2)$ :
&$i= 0,1,2$ \},\\

 $\A_1 = \{$
&$  (\infty_{1}, 3i+6, 3i+8)$,
&$  (\infty_{1}, 3i+6, 3i+5)$,
&$  (\infty_{1}, 3i+6, 3i+2)$,\\
&$  (\infty_{2}, 3i+2, 3i+6)$,
&$  (\infty_{2}, 3i+8, 3i+6)$,
&$  (\infty_{2}, 3i+5, 3i+6)$,\\
&$  (3i, 3i+7, 3i+8)$,
&$  (3i, 3i+1, 3i+5)$,
&$  (3i, 3i+4, 3i+2)$,\\
&$  (3i+6, 3i+3, 3i+7)$,
&$  (3i+6, 3i+3, 3i+1)$,
&$  (3i+6, 3i+3, 3i+4)$,\\
&$  (3i+7, 3i+2, 3i+5)$,
&$  (3i+1, 3i+8, 3i+2)$,
&$  (3i+4, 3i+5, 3i+8)$ :
&$i= 0,1,2 \}$,\\

 $\A_2 = \{$
&$  (\infty_{1}, 3i+3, 3i+4)$,
&$  (\infty_{1}, 3i+3, 3i+7)$,
&$  (\infty_{1}, 3i+3, 3i+1)$,\\

&$  (\infty_{2}, 3i+1, 3i+3)$,
&$  (\infty_{2}, 3i+4, 3i+3)$,
&$  (\infty_{2}, 3i+7, 3i+3)$,\\

&$  (3i, 3i+4, 3i+5)$,
&$  (3i, 3i+7, 3i+2)$,
&$  (3i, 3i+1, 3i+8)$,\\

&$  (3i+3, 3i+6, 3i+5)$,
&$  (3i+3, 3i+6, 3i+2)$,
&$  (3i+3, 3i+6, 3i+8)$,\\

&$  (3i+7, 3i+1, 3i+5)$,
&$  (3i+1, 3i+4, 3i+2)$,
&$  (3i+4, 3i+7, 3i+8)$ :
&$i= 0,1,2 \}$.\\
\end{tabular}}

{\small
\begin{tabular}{llllllllllllll}

 $A(D_0) =\{$
&(1,2), &(1,5), &(1,8), &(2,4), &(2,7), &(4,5), &(4,8), &(5,7), &(7,8), \\
&(2,1), &(4,2), &(5,1) ,&(5,4), &(7,2), &(7,5), &(8,1), &(8,4), &(8,7) $\}$,\\
$A(D_1) =\{$
&(0,2), &(0,5), &(0,8), &(2,3), &(2,6), &(3,5), &(3,8), &(5,6), &(6,8), \\
&(2,0), &(3,2), &(5,0), &(5,3), &(6,2), &(6,5), &(8,0), &(8,3), &(8,6) $\}$,\\
$A(D_2) =\{$
&(0,1), &(0,4), &(0,7), &(1,3), &(1,6), &(3,4), &(3,7), &(4,6), &(6,7), \\
&(1,0), &(3,1), &(4,0), &(4,3), &(6,1), &(6,4), &(7,0), &(7,3), &(7,6) $\}$.\\
\end{tabular}}

Let $G_r=\{\infty_{1},\infty_{2}, r,3+r,6+r\}$ and $\B_r= \{G_r,G_r, G_r\}\cup \A_r$ for all $r\in R$ .
It is easy to check that $(S, \infty_{1}, \infty_{2},\B_r, D_r)_{r\in R}$
is a GLS*$(1,3;2,(3, \{3\}, \{3\}, \{5\}), 11)$.\qed

The proofs of the following two recursive constructions for LS and GLS are similar to the proofs of Constructions 2.1 and 2.2$(c)$ in \cite{T1993}.

\begin{construction}\label{w=w}
Let $w$ be an odd positive integer.
 If a GLS$(1,\lambda; 2, (3, \{3\}, \{3\},\{3\}), v + 2)$ exists,
then there exists a GLS$(1,\lambda; 2,(3,\{3\},\{3\},\{w+2\}), wv +2)$.
\end{construction}

\begin{construction}\label{Dasui}
 If there exist a GLS*$(1,\lambda; 2, (3, K_0, \{3\}, K_2), v)$,
 an LS*$(1,\lambda; 2,(3,K_0^{\prime}), k)$ for all $k\in K_0$, and
 a GLS*$(1,\lambda; 2,(3,K_0^{\prime},\{3\}, K_2^{\prime} ), k)$ for all $k\in K_2$,
then there exists a GLS*$(1,\lambda; 2,$ $(3,K_0^{\prime},\{3\},K_2^{\prime}), v)$.
\end{construction}

\begin{lemma}\label{v=k}
1.  There exists an LS$(1,2;2,(3, \{3, 5\}),v)$ for $v=14,26$.

\noindent 2.  Let $w$ be an odd positive integer.
Then  there is an LS$(1,3;2,(3, \{3\},\{3\},\{w+2\}), 3w+2)$.
\end{lemma}

\noindent {\it Proof:} 1. By Lemmas \ref{v=6} and \ref{v=10} there is a GLS$(1,2;2,(3, \{3\}, \{3\}, \{3\}), \frac{v-2}{3}+2)$.
Apply Construction \ref{w=w} with $w=3$ to get an LS$(1,2;2,(3, \{3, 5\}),v)$ for $v=14,26$.

2. By Lemma \ref{v=5} there is a GLS$(1,3;2,(3, \{3\}, \{3\}, \{3\}), 5)$.
Apply Construction \ref{w=w} with odd positive integer $w$ to get an LS$(1,3;2,(3, \{3\},\{3\},\{w+2\}), 3w+2)$.\qed

\begin{theorem}\label{(3,2)}
1. There exists an LS*$(1,2;2,(3,\{3,8\}),v)$ for $v\equiv 8\pmod {12}$.

\noindent 2. There exists an LS*$(1,2;2,(3,\{3,14,26,50,98\}),v)$ for $v\equiv 2\pmod {12}$.
\end{theorem}

\noindent {\it Proof:} 1.  By Theorem~\ref{GLS} there is a GLS$(2,(3,\{3\},\{3\},\{5\}),\frac{ v+2}{2})$.
Apply Construction \ref{LS-LGDD}(3) with a $(3,1)$-LGDD$(2^3)$ to obtain an LS*$(1,2;2,(3,\{3,4,8\}),v)$.
An LS$(1,2;2,(3,\{3\}),4)$ exists by Theorem~\ref{LGDD}.
Then we apply Construction \ref{Dasui} to get an LS*$(1,2;2,(3,\{3,8\}),v)$.

2.  By Theorem~\ref{GLS} we have a GLS$(2,(3, \{3,4\}, \{3\}, \{8,14,26,50\}), \frac{ v+2}{2})$.
By Theorem~\ref{LGDD} there exist a $(3,1)$-LGDD$(2^3)$ and a $(3,1)$-LGDD$(2^4)$.
Apply Construction \ref{LS-LGDD} (3) to obtain an LS*$(1,2;2,(3,\{3,4,14,26,50,98\}),v)$.
Similarly,  we apply Construction \ref{Dasui} to get the required LS*$(1,2;2,(3,\{3,14,26,50,98\}),v)$.
\qed

Next we will give a recursive construction for LS*. For our construction, we need the following definitions.
 A {\it quasigroup} of order $n$ is a pair $(\mathbb{Q}, \circ)$,
 where $\mathbb{Q}$ is a set of size $n$ and $``\circ"$ is a binary operation on $\mathbb{Q}$
 such that for every pair of elements $a, b\in \mathbb{Q}$,
 the equations $a\circ x=b$ and $y\circ a=b$ have unique solutions.
 A quasigroup is said to be {\it idempotent} if $a\circ a=a$ for any $a\in  \mathbb{Q}$.
 A quasigroup is said to be {\it commutative} if $a\circ b=b\circ a$ for any $a, b \in  \mathbb{Q}$.
 An idempotent commutative quasigroup is a rename the table for $(\mathbb{Z}_{n}, +)$ which is the additive group of integers modulo $n$.
 For any odd $n$, there exists an idempotent commutative quasigroup  of order $n$.

\begin{construction}\label{w=SLS}
Let $w$ be an odd positive integer.
If there exists a GLS*$(1,\lambda; 2, (3, \{3\}, $ $\{3\},K_2), v + 2)$ with $k\ge 4$ for all $k\in K_2$,
then there is an LS*$(1,\lambda; 2,(3,\{3\},$ $\{3\},K_2^{\prime}), wv +2)$,
where $K_2^{\prime}=\{w(k-2)+2 : k\in K_{2}\}$.
\end{construction}
\noindent {\it Proof:}
Let  $(S, \infty_{1}, \infty_{2}, {\cal B}_r, D_r)_{r\in R}$ be a
GLS*$(1,\lambda; 2, (3, \{3\}, $ $\{3\},K_2), v + 2)$ and
$\{\infty_{1}, \infty_{2}\}\cap (S\times \mathbb{Z}_w)= \emptyset$.
Then $|R|=\frac{v}{\lambda}$.

{\bf Step 1:} For each $B\in {\cal B}_r$ and $\{\infty_{1}, \infty_{2}\}\subset B$,
let $\B_{B_{(r, i)}}= \{\{\infty_{1}, \infty_{2}\}\cup ((B\backslash\{\infty_{1}, \infty_{2}\})\times \mathbb{Z}_{w})\}$,  $i\in \mathbb{Z}_{w}$.

{\bf Step 2:}  For each $B\in {\cal B}_r$ and $\{\infty_{1}, \infty_{2}\}\cap B=\emptyset$,
let $(B\times \mathbb{Z}_w, \G_{B}, \B_{B_{(r,i)}})_{i\in \mathbb{Z}_{w}}$
be  a $(3,1)$-LGDD$(w^{3})$, where
$\G_{B}=\{\{x\}\times\mathbb{Z}_w : x\in B\}$.

{\bf Step 3:}
Let $(\mathbb{Z}_{w},\circ)$ be an idempotent commutative quasigroup  of order $w$.
For each $B=\{\infty_{l},x,y\}\in {\cal B}_r$, $l\in\{1, 2\}$, and  $(x,y)\in D_r$,
put ${\cal B}_{B_{(r,i,0)}}=\{\{\infty_{l},x_a,y_{a+i}\} : a \in\mathbb{Z}_{w} \}$ and
${\cal B}_{B_{(r,i,1)}}=\{\{x_{a},x_{b},y_{a \circ b +i}\} : a,b\in\mathbb{Z}_{w},a\neq b \}$.
Let  $\B_{B_{(r, i)}}={\cal B}_{B_{(r,i,0)}}\cup {\cal B}_{B_{(r,i,1)}}$,
$ i\in \mathbb{Z}_{w}$.

Let  $\B_{(r,i)}= \cup_{B\in {\cal B}_r} \B_{B_{(r, i)}}$, $(r, i)\in R_1=R\times \mathbb{Z}_{w}$.
We shall show $( S\times \mathbb{Z}_w,\infty_{1}, \infty_{2}, \B_{(r,i)})_{(r,i)\in R_1}$
is the desired LS*$(1,\lambda; 2,(3,\{3\},\{3\},K_2^{\prime}), wv +2)$.

Firstly, we show $( S\times \mathbb{Z}_w,\infty_{1}, \infty_{2}, \B_{(r,i)})$
is an S$(\lambda; 2,(3,\{3\},\{3\},K_2^{\prime}), wv +2)$.
Let $P$ be a 2-subset of $(S\times \mathbb{Z}_w)\cup \{\infty_{1}, \infty_{2}\}$.
We distinguish 4 cases.

(1) $P=\{\infty_{1}, \infty_{2}\}$.
There are  exactly $\lambda$ identical blocks $B\in {\cal B}_r$ such that $P\subset B$ since $(S, \infty_{1}, \infty_{2}, {\cal B}_r, D_r)_{r\in R}$ is a
GLS*. By Step 1 $P$ appears in $\lambda$ identical blocks $\{\infty_{1}, \infty_{2}\}\cup ((B\backslash\{\infty_{1}, \infty_{2}\})\times \mathbb{Z}_{w})$ in $\B_{(r,i)}$.

(2)  $P=\{\infty_{l}, x_a\}, l\in\{1, 2\},  x \in S$ and $j\in \mathbb{Z}_w$.
Let multiset $M_1=\{B: B\in {\cal B}_r, \{\infty_{l}, x\} \subset B \}$.
Then we have $|M_1|=\lambda$.
If there is a block $B\in M_1$ such that $\{\infty_{1}, \infty_{2}\}\subset B$,
then $|B|\ge 4$ since $k\ge 4$ for all $k\in K_2$, and  $M_1$ is consisted of $\lambda$ identical blocks by the definition of a GLS*.
By Step 1 $P$ appears in $\lambda$ identical blocks $\{\infty_{1}, \infty_{2}\}\cup ((B\backslash\{\infty_{1}, \infty_{2}\})\times \mathbb{Z}_{w})$ in $\B_{(r,i)}$.
Otherwise, for any block $B\in M_1$, we have $\{\infty_{1}, \infty_{2}\}\not\subset B$.
Suppose $B=\{\infty_{l}, x,y\}$.
Then by Step 3 $P$ is contained in the block $\{\infty_{l}, x_a,y_{a+i}\}\in {\cal B}_{B_{(r,i,0)}}$ if $(x,y)\in D_r$ or  in the block $\{\infty_{l}, y_{a-i},x_{a}\}\in {\cal B}_{B_{(r,i,0)}}$  if $(y,x)\in D_r$.
So $P$ is contained in exactly $\lambda$ blocks of $\B_{(r,i)}$ since $|M_1|=\lambda$.

(3) $P=\{x_a, x_b\}, x \in S$, $a, b\in \mathbb{Z}_w$ and $a\neq b$.
If there is a block $B\in {\cal B}_r$ such that $\{\infty_{1}, \infty_{2},x\}\subset B$, then $t_x=\lambda$ and by Step 1 $P$ appears in $\lambda$ identical blocks in $\B_{(r,i)}$.
Otherwise, $t_x=0$. Let $M_2= \{y: (x,y)\in D_r \}$.
Then $|M_2|=d^+(x)=\lambda-t_x=\lambda$.
For any $y \in M_2$, there is a block $B=\{\infty_{l}, x,y\}\in {\cal B}_r$, $l\in\{1,2\}$.
Then by Step 3 $P$ is contained in the block $\{x_a, x_b, y_{a\circ b+i}\} \in {\cal B}_{B_{(r,i,1)}}$.
So $P$ is contained in exactly $\lambda$ blocks of $\B_{(r,i)}$ since $|M_2|=\lambda$.

(4)  $P=\{x_a, y_b\}, x, y \in S$, $x\neq y$ and $a,b\in \mathbb{Z}_w$.
Let multiset $M_3=\{B: B\in {\cal B}_r, \{x, y\} \subset B \}$.
Then $|M_3|=\lambda$.
If there is a block $B\in M_3$ and $|\{\infty_{1}, \infty_{2}\}\cap B|=2$,
then $P$ appears in $\lambda$ identical blocks in $\B_{(r,i)}$ from Step 1.
Otherwise, for any block $B\in M_3$, $|\{\infty_{1}, \infty_{2}\}\cap B|\le 1$.
We distinguish into 2 subcases.

(I) $|\{\infty_{1}, \infty_{2}\}\cap B|= 1$.
Let $B=\{\infty_{l}, x,y\}$, $l\in \{1,2\}$ and $(x,y)\in D_r$.
If $i= b-a$, then $P$ is contained in the block $\{\infty_{l}, x_a,y_{b}\}\in {\cal B}_{B_{(r,i,0)}}$ from Step 3.
If $i\neq b-a$, then there is a unique $c\in\mathbb{Z}_w$ such that $b-i=a\circ c$ by idempotent commutative quasigroup $(\mathbb{Z}_{w},\circ)$.
Then $P$ is contained in the block $\{x_a, x_c, y_b\}\in {\cal B}_{B_{(r,i,1)}}$ from Step 3.

(II) $|\{\infty_{1}, \infty_{2}\}\cap B|= 0$.
 Suppose  $B=\{ x,y, z\}$.
By Step 2  $P$ is contained in exactly one block of  $\B_{B_{(r,i)}}$ since  $(B\times \mathbb{Z}_w, \G_{B}, \B_{B_{(r,i)}})$
is a $(3,1)$-GDD$(w^{3})$.
So $P$ is contained in exactly $\lambda$ blocks of $\B_{(r,i)}$ since $|M_3|=\lambda$.

Next we prove that for any 3-subset $T$ of $(S\times \mathbb{Z}_w)\cup \{\infty_{1}, \infty_{2}\}$, there is exactly one block $A\in \cup_{(r,i)\in R} \B_{(r,i)}$ such that $T\subseteq A$.
We distinguish 6 cases.

(1) $T=\{\infty_{1}, \infty_{2}, x_a\}, x \in S$ and $a\in \mathbb{Z}_w$.
Then $\{\infty_{1}, \infty_{2}, x\}$ is contained in a unique block $B \in \cup_{r\in R}{\cal B}_r$  since $(S, \infty_{1}, \infty_{2}, {\cal B}_r, D_r)_{r\in R}$ is a GLS*.
By Step 1 $T$ is contained in the block $A = \{\infty_{1}, \infty_{2}\}\cup ((B\backslash\{\infty_{1}, \infty_{2}\})\times \mathbb{Z}_{w})$
in $\cup_{(r,i)\in R_1} \B_{(r,i)}$.

(2)  $T=\{\infty_{l}, x_a, x_b\}, l\in\{1, 2\},  x \in S, a, b\in \mathbb{Z}_w$ and $a\neq b$. Then this case is similar to case (1).

(3)  $T=\{\infty_{l}, x_a, y_b\}, l\in\{1, 2\},  x,y \in S$, $x\neq y$ and $a, b\in \mathbb{Z}_w$.
Then there is a unique block $B\in \cup_{r\in R}{\cal B}_r$ such that $\{\infty_{l}, x, y\}\subset B$.
 If $|\{\infty_{1},\infty_{2}\}\cap B| =2$,
then this case is similar to case (1).
Otherwise,  $|\{\infty_{1},\infty_{2}\}\cap B| =1$.
Then there is a unique $r\in R$ such that $ B=\{\infty_{l}, x, y\}\in {\cal B}_r$.
If $(x,y)\in D_r$, then there is a unique $i_1= b-a$ such that
 $T \in {\cal B}_{B_{(r,i_1,0)}}$ from Step 3.
If $(y,x)\in D_r$, then there is a unique $i_1= a-b$ such that
 $T \in {\cal B}_{B_{(r,i_1,0)}}$ from Step 3.
So $T$ is contained in exactly one blocks of $\cup_{(r,i)\in R_1} \B_{(r,i)}$.

(4) $T=\{x_a, x_b, x_c\}, x \in S$, $a, b, c\in \mathbb{Z}_w$ and $|\{a, b, c\}| =3$.
This case is similar to case (1).

(5) $T=\{x_a, x_b, y_c\}, x, y \in S$, $x\neq y$, $a, b, c\in \mathbb{Z}_w$ and $a\neq b$.
 If there is a unique block $B$ such that $\{\infty_{1}, \infty_{2}, x,y\}\subset B$ and $B\in \cup_{r\in R}{\cal B}_r$,
then this case is similar to case (1).
Otherwise, there is a unique block $B=\{\infty_{l}, x,y\}$, $l\in \{1,2\}$ such that $(x,y)\in \cup_{r\in R}D_r$.
Then there is a unique  $i =c - a\circ b$ such that $T \in {\cal B}_{B_{(r,i,1)}}$ from Step 3.
So $T$ is contained in exactly one blocks of $\cup_{(r,i)\in R_1} \B_{(r,i)}$.

(6)  $T=\{x_a, y_b, z_c\}, x, y, z \in S$, $|\{x, y, z\}| =3$ and $a,b,c\in \mathbb{Z}_w$.
Then there is a unique block $B\in {\cal B}_r$ such that $\{x, y, z\}\subset B$.
If $|\{\infty_{1},\infty_{2}\}\cap B| =2$,
then this case is similar to case (1).
Otherwise, $|\{\infty_{1},\infty_{2}\}\cap B| \le 1$. Then we have $|\{\infty_{1},\infty_{2}\}\cap B| =0$.
By Step 2 we have $T\in \cup_{i\in \mathbb{Z}_{w}}\B_{B_{(r,i)}}$ since  $(B\times \mathbb{Z}_w, \G_{B}, \B_{B_{(r,i)}})_{i\in \mathbb{Z}_{w}}$
is a $(3,1)$-LGDD$(w^{3})$.
So $T$ is contained in exactly one blocks of $\cup_{(r,i)\in R_1} \B_{(r,i)}$.

Finally, we show that each block $A$,  $A\in\B_{(r,i)}$ and $|A|\ge 4$, appears $\lambda$ times in $\B_{(r,i)}$ and $|A|-2$ times in the multiset $\{\B_{(r,i)}:(r,i)\in R_1\}$, respectively.

Let $A\in\B_{(r,i)}$ and $|A|\ge 4$.
Then $A$ must come from Step 1.
So we have  $\{\infty_{1}, \infty_{2}\} \subset A$.
Thus, there is a block  $B\in \B_r$ such that $\{\infty_{1}, \infty_{2}\} \subset B$ and $A= \{\infty_{1}, \infty_{2}\}\cup ((B\backslash\{\infty_{1}, \infty_{2}\})\times \mathbb{Z}_{w})$ from Step 1.
 Since $(S, \infty_{1}, \infty_{2}, {\cal B}_r, D_r)_{r\in R_1}$ is a
GLS*,  $B$ appears $\lambda$ times in $\B_{r}$ and $|B|-2$ times in the multiset $\{\B_{r}:r\in R\}$, respectively.
So  $A$ appears $\lambda$ times in $\B_{(r,i)}$ and $w(|B|-2)=|A|-2$ times in the multiset $\{\B_{(r,i)}:(r,i)\in R_1\}$ from Step 1, respectively.

Now the proof is complete. \qed

\section{Constructions for GLS*}

In this section, we shall give two constructions for good large sets by using LR designs and generalized frames.

A GDD $(X,{\cal G},{\cal B})$ is called {\it resolvable} if there exists a
partition $\Gamma = \{P_1, P_2,\ldots, P_r \}$ of ${\cal B}$ such that each part $P_i$ (called a {\it parallel class}) is itself a partition of $X$.
The  partition $\Gamma$ is called a {\it resolution} of ${\cal B}$.

Let $X$ be a $v$-set. An LR {\it design} of order $v$ (or LR$(v)$ as in \cite{L2002}) is a collection
$\{(X, \A_{k}^{j}) : 1\le k\le \frac{v-1}{2},  j=0,1\}$ of S$(2,3,v)$s with following properties:
(1) Let the resolution of $\A_{k}^{j}$ be $\Gamma_{k}^{j} = \{A_{k}^{j}(h) : 1\le h\le \frac{v-1}{2}\}$.
There is an element in each $\Gamma_{k}^{j}$,
 which without loss of generality, we can suppose is $A_{k}^{j}(1)$,
such that
$$\bigcup_{k=1}^{\frac{v-1}{2}} A_{k}^{0}(1) = \bigcup_{k=1}^{\frac{v-1}{2}}A_{k}^{1}(1)= \A,$$
and $(X, \A)$ is an S$(2,3,v)$.
(2) For any triple $T = \{x,y,z\}\subset X$,  $|\{x,  y, z\}|=3$,
there exist $k, j$ such that $T\in \A_{k}^{j}$.

\begin{theorem}\label{LR}{\rm (\cite{JL2008,L2002})}
There exists an LR$(v)$ for any $v\in\{3^n,2 \times 7^n + 1,2 \times 13^n + 1:n\ge 1\}$.
\end{theorem}

\begin{construction}\label{CLR}
If there exists an LR$(2v+1)$,
then there exist a GLS$( 2, (3, \{3, 6\}, \{3\}, $ $\{6\}), 4v+2)$
and a GLS*$( 1,2;2, (3, \{3, 6\}, \{3\}, \{6\}), 4v+2)$.
\end{construction}

\noindent {\it Proof:}
Let $S$ be a $2v$-set and $\infty \notin S$.
Let $\{(S\cup\{\infty\}, \A_{k}^{j}) : 1\le k\le v,  j=0,1\}$ be a  LR$(2v+1)$,
where each $\A_{k}^{j}$ can be partitioned into parallel classes
$\{A_{k}^{j}(h) : 1\le h\le v\}$, with
$$\bigcup_{k=1}^{v} A_{k}^{0}(1) = \bigcup_{k=1}^{v}A_{k}^{1}(1)= \A,$$
such that $(S\cup\{\infty\}, \A)$ is an S$(2,3,2v+1)$.

{\bf Step 1:} For each $B\in A_{k}^{j}(1)$,
let $\B_{B_i}= \{ B\times \mathbb{Z}_2\}, i\in\{0,1\}$.

{\bf Step 2:} For each  $B\in A_{k}^{j}\backslash A_{k}^{j}(1)$ and $\infty \notin B$,
let $(B\times \mathbb{Z}_2,\{\{x\}\times \mathbb{Z}_2 : x \in B\}, \B_{B_i})_{i\in\{0,1\}}$
be a $(3,1)$-LGDD$(2^3)$.

{\bf Step 3:} For each $B\in A_{k}^{j}\backslash A_{k}^{j}(1)$ and $\infty \in B$,
let $B=\{\infty, x, y\}$ and

$\B_{B_0}=\{(\infty_1,(x,0),(y,0)), (\infty_1,(x,1),(y,1)), (\infty_2,(y,1),(x,0)), (\infty_2,(y,0),(x,1))\}$,

$\B_{B_1} = \{(\infty_1,(x,0),(y,1)), (\infty_1,(x,1),(y,0)), (\infty_2,(y,0),(x,0)), (\infty_2,(y,1),(x,1))\}$.

\noindent It is easy to check that  $( B\times \mathbb{Z}_2, \{\{x\}\times \mathbb{Z}_2 : x \in B\}, \B_{B_i})_{i\in\{0,1\}}$
is a $(3,1)$-LGDD$(2^3)$.

Let $V(D_{B_i})=\{x,y\}\times \mathbb{Z}_2$ and

$A(D_{B_0}) = \{((x,0),(y,0)), ((x,1),(y,1)), ((y,1),(x,0)), ((y,0),(x,1))\}$,

$A(D_{B_1}) = \{((x,0),(y,1)), ((x,1),(y,0)), ((y,0),(x,0)), ((y,1),(x,1))\}$.

\noindent Then we have that $D_{B_i}=(V(D_{B_i}), A(D_{B_i}))$, $i=0,1$, are two Eulerian digraphs.

Let $R = \{0,1\}\times \{0,1\} \times\{1,2,\ldots, v\}$.
For each $(i,j,k)\in R$,  construct
 a digraph $D_{ijk}$ on $S\times \mathbb{Z}_2$ such that
$A(D_{ijk}) = \{A(D_{B_{i}}) : B\in \A_{k}^{j} \setminus \A_{k}^{j}(1),  \infty \in B \}$.
Then  $D_{ijk}$  is an Eulerian digraph since $D_{B_{i}}$ is an Eulerian digraph.
Let $\B_{ijk} = \cup_{B\in \A_{k}^{j}}\B_{B_{i}}$.
We shall show $( S\times \mathbb{Z}_w,\infty_{1}, \infty_{2}, \B_{ijk},D_{ijk})_{(i,j,k)\in R}$
is the desired GLS$( 2,(3,\{3,6\},\{3\},\{6\}), 4v+2)$.

Firstly, we shall show $(S\times \mathbb{Z}_2,\infty_1,\infty_2,\B_{ijk}, D_{ijk})$
 is a GS$(2,(3, \{3,6\}, \{3\}, \{6\}), 4v+2)$.

Let $P=\{(x,a),(y,b)\}$ be a 2-subset of $(S\times \mathbb{Z}_2) \cup \{\infty_1,\infty_2\}$.
If $x=y$,   then there is exactly one block $B\in A_{k}^{j}(1)$ such that $x\in B$ since $A_{k}^{j}(1)$ is a parallel class.
 By Step 1 $P$ is contained in the block $B\times \mathbb{Z}_2$.
If $x\neq y$.  There is exactly one block $B\in A_{k}^{j}$ such that $\{x,y\}\subset B$ since  $(S\cup\{\infty\}, A_{k}^{j})$ is an S$(2,3,2v+1)$.
 By Steps 2 and 3  $P$ is contained in exactly one block of  $\B_{B_{i}}$ since  $(B\times \mathbb{Z}_2,\{\{x\}\times \mathbb{Z}_2 : x \in B\}, \B_{B_i})_{i\in\{0,1\}}$ is a $(3,1)$-LGDD$(2^3)$.

For any $(x,a)\in S\times \mathbb{Z}_2$, let $t_{(x,a)}=|\{A : \{\infty_{1}, \infty_{2}, (x,a)\} \subset A, A \in{\cal B}_{ijk} \}|$ and
 $d^+(x,a)=|\{((x,a),(y,b)) : (y,b)\in S\times \mathbb{Z}_2, ((x,a),(y,b)) \in D_{ijk}, \{\infty_{l}, (x,a),(y,b)\}$ $ \in {\cal B}_{ijk}, l= 1,2\}|$.
 Since $D_{ijk}$ is an Eulerian digraph, we have $d^-(x,a)=d^+(x,a)$.
There is a unique block $B\in A_k^j$ such that $\{x,\infty\} \subset B$ since $(S\cup \{\infty\}, A_k^j)$ is an S$(2,3, 2v+1)$.
If  $B\in A_k^j(1)$, we have $\{\infty_{1}, \infty_{2}, (x,a)\} \subset B\times \mathbb{Z}_2$ from Step 1.
Then $t_{(x,a)}=1 $ and $d^-(x,a)=d^+(x,a)=0$.
If $B\in A_k^j\backslash A_k^j(1)$, there is a unique point $(y,b)\in S\times \mathbb{Z}_2$ such that $((x,a),(y,b)) \in A(D_ijk)$  and $\{\infty_{l}, (x,a), (y,b)\}\in \B_{B_i}$, $l\in \{1,2\}$  from Step 3.
Then $t_{(x,a)}=0$ and $d^-(x,a)=d^+(x,a)=1$.

Next we prove that for any 3-subset $T=\{(x,a),(y,b),(z,c)\}$ of $(S\times \mathbb{Z}_2)\cup \{\infty_{1}, \infty_{2}\}$, there is exactly one block $A\in \cup_{(i,j,k)\in R} \B_{ijk}$ such that $T\subseteq A$.
We distinguish 2 cases.

(1) $|\{x,y,z\}|=2$. Then we suppose that $x\neq y=z$.
There is exactly one block $B\in \A$ such that $\{x,y\}\subset B$ since $(S\cup\{\infty\}, \A)$ is an S$(2,3,2v+1)$.
By Step 1 $T$ is contained in the block $B\times \mathbb{Z}_2$.

(2) $|\{x,y,z\}|=3$. By definition of LR design, there exist $k, j$ such that $B=\{x,y,z\}\in \A_{k}^{j}$.
If $B\in A_{k}^{j}(1)$, this case is similar to case (1).
Otherwise, $B\in A_{k}^{j}\backslash A_{k}^{j}(1)$.
By Steps 2 and 3 $T  \in \B_{B_0}\cup \B_{B_1}$.

Thirdly, we show that each block $A$,  $A\in\cup_{(i,j,k)\in R}\B_{ijk}$ and $|A|=6$,
appears $4$ times in the multiset $\{\B_{ijk}:(i,j,k)\in R\}$.

Let $A\in\cup_{(i,j,k)\in R}\B_{ijk}$ and $|A|=6$. Then $A$ must come from Step 1.
So there is a unique block $ B\in\A$ such that $A=B\times\mathbb{Z}_2$.
By Step 1  $B$ appears exactly $4$ times in the multiset $\{\B_{ijk}:(i,j,k)\in R\}$ since
  $\cup_{k=1}^{v} A_{k}^{0}(1) = \cup_{k=1}^{v} A_{k}^{1}(1)= \A$.

Finally, we show that  each ordered pair $((x,a),(y,b))$ of distinct elements of $S\times \mathbb{Z}_2$,
not contained in some block $A$ with $\{\infty_{1},\infty_{2}\}\subset A$, appears in exactly one $D_{ijk}$.

Let ordered pair $((x,a),(y,b))$ of distinct elements of $S\times \mathbb{Z}_2$,
not contained in some block $A$ with $\{\infty_{1},\infty_{2}\}\subset A$.
Then we have $x\neq y$.
There exist $k, j$ such that $\{x,y,\infty\}\in \A_{k}^{j}$.
By Step 3 we have $((x,a),(y,b))\in D_{ojk} \cup D_{1jk}$.

Thus, $( S\times \mathbb{Z}_2,\infty_{1}, \infty_{2}, \B_{ijk},D_{ijk})_{(i,j,k)\in R}$
is the desired GLS$( 2,(3,\{3,6\},\{3\},\{6\}), 4v+2)$.

Let $\B_{jk}= \B_{0jk} \cup \B_{1jk}$, $D_{jk}= D_{0jk} \cup D_{1jk}$ and $R_1=\{0,1\}\times\{1,2,\ldots,v\}$.
It is easy to check that $(S\times \mathbb{Z}_2,\infty_1,\infty_2,\B_{jk}, D_{jk})_{(j,k)\in R_1}$
 is a GLS*$(1,2; 2,(3,\{3,6\},\{3\},\{6\}),4v+2)$.

Now the proof is complete. \qed

By Theorem \ref{LR} and Construction \ref{CLR}, we have the following theorem.
\begin{theorem}\label{w=LR}
There exist a GLS$( 2, (3, \{3, 6\}, \{3\}, \{6\}), 2v)$ and
a GLS*$(1,2; 2, $ $(3,\{3, 6\},\{3\},$ $ \{6\}), 2v)$ for all $v\in\{3^n,2 \times 7^n + 1,2 \times 13^n + 1:n\ge 1\}$.
\end{theorem}

\begin{lemma}\label{v=50}
 There exists an  LS*$(1,2; 2, (3, \{3\}, \{3\}, \{14\}), 50)$.
\end{lemma}

\noindent {\it Proof:}
By Theorem \ref{w=LR} there is a GLS*$(1,2; 2, (3, \{3,6\}, \{3\}, \{6\}), 18)$.
Apply Construction \ref{Dasui} with an LS$(1,2;2,(3, \{3\}), 6)$ from Lemma \ref{v=6} to
get a GLS*$(1,2; 2, (3, \{3\}, \{3\}, \{6\}), 18)$.
Then  we apply  Construction \ref{w=SLS} with $w=3$ to get
 an LS*$(1,2; 2,(3,\{3\}, \{3\}, \{14\}), 50)$.\qed

\begin{lemma}\label{v=98}
 There exists an LS*$(1,2; 2, (3, \{3\}, \{3\}, \{26\}),  98)$.
\end{lemma}
\noindent {\it Proof:}
By Theorem \ref{w=LR} there is a GLS$(2, (3, \{3,6\}, \{3\}, \{6\}), 18)$ $(S,\infty_{1}, \infty_{2},\beta_r,D_r)_{r\in R}$.

For $B\in \beta_r$ and $\{\infty_{1}, \infty_{2}\}\subset B$,
put $\B_{B_{(r, i)}}= \{\{\infty_{1}, \infty_{2}\}\cup ((B\backslash\{\infty_{1}, \infty_{2}\})\times \mathbb{Z}_{6})\}$,
$ i\in \mathbb{Z}_6$.

 For $B\in \beta_r$ and $\{\infty_{1}, \infty_{2}\}\cap B=\emptyset$,
let $R_B =\{\{r\}\times \mathbb{Z}_6 : B\in \beta_r \}$ .
By Theorem~\ref{LGDD}, we have a $(3,1)$-LGDD$(6^3)$ and $(3,1)$-LGDD$(6^6)$.
Then let $(B\times \mathbb{Z}_6, \G_{B}, \B_{B_{(r,i)}})_{(r,i)\in R_{B}}$
be  a $(3,1)$-LGDD$(6^{|B|})$, where
$\G_{B}=\{\{x\}\times\mathbb{Z}_6 : x\in B\}$.

Let
$A_0=\{(0,1,2),(0,4,5),(3,1,5),(3,4,2),(1,0,2),(1,3,5),(4,0,5), $ $(4,3,2),$ $ (2,0,1),$ $(2,3,4),(5,0,4),$ $(5,3,1)\}$,
$A_1=\{(0,0,2),(0,3,5),(3,0,5),(3,3,2),(1,0,1),(1,3,4),(4,0,4),$ $(4,3,1),(2,1,2),$ $(2,4,5),(5,1,5),(5,4,2)\}$,
$A_2=\{(0,0,1),(0,3,4),(3,0,4),(3,3,1),$ $(1,1,2),$ $(1,4,5),(4,1,5),$ $(4,4,2),(2,0,2),$ $(2,3,5),(5,0,5),(5,3,2)\}$.

For $B\in \beta_r$,  $B=\{\infty_{1},x,y\}$ and $(x,y)\in D_r$,
let $\B_{B_{(r,i)}}$, $0\le i\le 2$, consist of the sets
$\{\{(x,a),(x,b),(y,c)\},\{(x,a),(x,b),(y,c+3)\}\}, (a,b,c)\in A_i$,
as well as the sets
$\{\{(x,a),(x,a+3),(y,a+i),(y,a+i+3)\},\{(x,a),(x,a+3),(y,a+i),(y,a+i+3)\} \}, a\in \mathbb{Z}_6$.

For $B\in \beta_r$,  $B=\{\infty_{2},x,y\}$ and $(x,y)\in D_r$,
let $\B_{B_{(r,i)}}$, $0\le i\le 2$, consist of the sets
$\{\{(x,a),(x,b),(y,c)\},\{(x,a),(x,b),(y,c+3)\}\}, (a,b,c)\in A_i$, as well as the sets
$\{\{\infty_{l}, $ $(x,a),(y,a+i)\}, \{\infty_{1},(x,a),(y,a+i+3)\}\}$, $l=1,2$.

Let $\B_{(r,i)}= \cup_{B\in \beta_r} \B_{B_{(r, i)}}$ for $r\in R$ and  $0\le i\le 2$.
It is easy to check that  $(S\times \mathbb{Z}_6,\infty_{1},\infty_{2},\beta_{(r,i)})_{(r,i)\in R\times \{0,1,2\}}$ is
an LS*$(1,2; 2,(3, \{3, 4\}, \{3\}, \{26\}), 98)$.
Apply  Construction \ref{Dasui} to get an LS*$(1,2; 2, (3, \{3\}, \{3\}, \{26\}),  98)$.\qed

Let $X$ is a $v$-set.
The $(s+3)$-tuple $(X,\G,\A_1,\ldots,\A_{s},{\cal T})$ is called an $s$-{\it fan design}
(as in \cite{H1994}) if $(X, \G)$ is a 1-design,  $(X, \G\cup \A_i)$ is a PBD for $1\le i\le s$
and $(X, (\G\cup {\cal T})\cup (\bigcup_{i=1}^{s}\A_i))$ is a 3-wise balanced  design.
If block sizes of $\A_i$ and ${\cal T}$ are from $K_i$ $(1\le i\le s)$ and $K_T$, respectively,
then the $s$-fan design is denoted by $s$-FG$(3,(K_1,\ldots,K_{s},K_{T}), v)$.

A {\it generalized frame} (as in \cite{T1994}) F$(3, 3, g^u)$ is a GDD$(3,3;g^u)$ $(X,\G,\A)$
 such that the block set $\A$ can be partitioned into $gn$ subsets $\A_ x$, $x\in G$ and $G\in \G$,
each $(X\backslash G, \G\backslash \{G\}, \A_x)$ being a $(3,1)$-GDD$(g^{u-1})$.

\begin{theorem}\label{GF} {\rm(\cite{J2009,M1990,T1994})}
For $u\ge 3$ and $u\neq 5$,  an F$(3, 3, g^u)$  exists
if and only if $gu$ is even and $g(u-1)(u-2)$ is divisible by 3.
 For $u=5$, an F$(3, 3, g^5)$ exists if $g$ is even, $g\neq2$, and $g \not\equiv10,26 \pmod {48}$.
\end{theorem}

\begin{construction}\label{PCS}
Suppose that there exists a 1-FG$(3,(K_1,K_{T}), u)$.
Suppose that there exists a  GLS*$(1,m;2,(3,K_0,\{3\},\{m+2\}),mk+2)$ for any $k \in K_1$,
and there exists an  F$(3, 3, m^k)$ for any $k\in K_T$.
Then there is a GLS*$(1,m;2,(3,\{3\}\cup K_0,\{3\},\{m+2\}),mu+2)$.
\end{construction}

\noindent {\it Proof:}
Let $(X,\G,\A_1,{\cal T})$ be the given 1-FG$(3,(K_1,K_{T}), u)$.
Let $S = \{\infty_1, \infty_2\}$ and $(X\times \mathbb{Z}_m)\cap S= \emptyset$.
We shall construct the desired design on
$( X\times \mathbb{Z}_m)\cup S$.
Denote $G_x =\{x\}\times\mathbb{Z}_m$ for $x\in X$.
Let $M_x=\{G_x\cup S, \ldots, G_x\cup S\}$ be an $m$-multiset.

{\bf Step 1:} For each block $B\in \A_1$,
let  $(B\times\mathbb{Z}_m,\infty_{1}, \infty_{2}, \B_{B_x}, D_{B_x} )_{x\in B}$
be a GLS*$(1,m;2,(3, $ $K_0,\{3\},\{m+2\}),m|B|+2)$ such that $M_x$ is a set of all blocks of size $m+2$ in $\B_{B_x}$.

{\bf Step 2:} For each block $B\in {\cal T}$,
we can construct a generalized frame F$(3,3,m^{|B|})$ on $B\times \mathbb{Z}_m$
having $\Gamma_{B} = \{G_x: x\in B\}$ as its group set and
 the block set can be partitioned into $m|B|$ disjoint block sets
$\C_B(x,i)$ $(x\in B, i\in \mathbb{Z}_m)$
with the property that each
$((B\backslash \{x\})\times \mathbb{Z}_m, \Gamma_{B}\backslash\{G_{x}\}, \C_B(x,i))$
  is a $(3,1)$-GDD$(m^{|B|-1})$. Let $\C_{B_x} = \cup_{i=0}^{m-1}\C_{B}(x,i)$.

{\bf Step 3:} For any $x\in X$, let
$$\F_x=M_x\cup (\bigcup_{x\in B, B\in \A_1}(\B_{B_x}\backslash M_x))\cup(\bigcup_{x\in B, B\in {\cal T}}\C_{B_x}),   ~~~~~
D_x = \bigcup_{x\in B,B\in \A_1} D_{B_{x}}.$$

We shall show that
 $(X\times\mathbb{Z}_m, \infty_{1}, \infty_{2}, \F_x, D_x )_{x\in X}$
is a GLS*$(1,m;2,(3,\{3\}\cup K_0,\{3\}, $ $\{m+2\}),mu+2)$.

Firstly, we shall show  $(X\times\mathbb{Z}_m, \infty_{1}, \infty_{2}, \F_x, D_x )$
is a GS$(1,m;2,(3,\{3\}\cup K_0,\{3\}, $ $\{m+2\}),mu+2)$.

Let $P=\{(y,a),(z,b)\}$ be a 2-subset of $(X\times\mathbb{Z}_m) \cup \{\infty_1,\infty_2\}$. We distinguish two cases.

(1) $y=z$.  If $y \not\in \{x,\infty\}$, then there is exactly one block $B\in \A_1$ such that $\{x,y\}\subset B$ since $(X,\A_1)$ is an S$(2,K_1,u)$.
 By Step 1 $P$ is exactly contained in $m$ blocks of $\B_{B_x}$ since
$(B\times\mathbb{Z}_m,\infty_{1}, \infty_{2}, \B_{B_x}, D_{B_x} )$
is a GS$(1,m;2,(3, $ $K_0,\{3\},\{m+2\}),m|B|+2)$.
Otherwise, $y \in \{x,\infty\}$. For any block  $B\in \A_1$ satisfying $x\in B$,  by Step 1 there are  exactly $m$ identical blocks $A=(G_x\cup S)\in \B_{B_x}$ since $(B\times\mathbb{Z}_m,\infty_{1}, \infty_{2}, \B_{B_x}, D_{B_x})$ is a GS$(1,m;2,(3, $ $K_0,\{3\},\{m+2\}),m|B|+2)$ and $M_x\subset \B_{B_x}$.

(2) $y\neq z$.  Then there is exactly one block $B\in \A_1\cup {\cal T}$ such that $\{x,y,z\}\subset B$
since  $(X,\G,\A_1,{\cal T})$ is a 1-FG$(3,(K_1,K_{T}), u)$.
If $B\in \A_1$, by Step 1 $P$ is exactly contained in $m$ blocks of $\B_{B_x}$ since
$(B\times\mathbb{Z}_m,\infty_{1}, \infty_{2}, \B_{B_x}, D_{B_x} )$ is a GS.
Otherwise,   $B\in {\cal T}$, by Step 2 $P$ is exactly contained in one block of $\C_B(x,i)$ since  $((B\backslash \{x\})\times \mathbb{Z}_m, \Gamma_{B}\backslash\{G_{x}\}, \C_B(x,i))$ is a $(3,1)$-GDD$(m^{|B|-1})$.
So $P$ is exactly contained in $m$ blocks of $\B_{B_x}$ since  $\C_{B_x} = \cup_{i=0}^{m-1}\C_{B}(x,i)$.

For any $(y,a)\in X\times\mathbb{Z}_m$, let $t_{(y,a)}=|\{A : \{\infty_{1}, \infty_{2}, (y,a)\} \subset A, A \in{\cal B}_{x} \}|$ and
 $d^+(y,a)=|\{((y,a),(z,b)) : (z,b)\in X\times\mathbb{Z}_m, ((y,a),(z,b)) \in D_{x}, \{\infty_{l}, (y,a),(z,b)\}$ $ \in {\cal B}_{x}, l= 1,2\}|$.
 Since $D_{B_x}$ is an Eulerian digraph, we have $D_{x}$ is an Eulerian digraph. So we have $d^-(y,a)=d^+(y,a)$.
 If $y=x$, we have $t_{(y,a)}=|M_x|=m$ and $d^-(x,a)=d^+(x,a)=0$ by Step 1.
 Otherwise $y\neq x$. Then there is a unique block $B\in \A_1$ such that $\{x,y\}\subset B$ since $(X,\A_1)$ is an S$(2,K_1,u)$.
 By Step 1, we have $d^+(y,a)=m$ and $t_{(x,a)}=0$ since $(B\times\mathbb{Z}_m,\infty_{1}, \infty_{2}, \B_{B_x}, D_{B_x} )_{x\in B}$ is a GS$(1,m;2,(3, $ $K_0,\{3\},\{m+2\}),m|B|+2)$.

Next we prove that for any 3-subset $T=\{(y,a),(z,b),(g,c)\}$ of $(X\times \mathbb{Z}_m)\cup \{\infty_{1}, \infty_{2}\}$, there is exactly one block $A\in \cup_{x\in X} {\cal F}_{x}$ such that $T\subseteq A$.
We distinguish 2 cases.

(1) $|T\cap (G_x\cup S)|=3$ for some $x\in X$.
For any $B\in \A_1$ and $x\in B$, by Step 1 $T$ is exactly contained in the block $G_x\cup S\in \cup_{x\in B}\B_x$ since
$(B\times\mathbb{Z}_m,\infty_{1}, \infty_{2}, \B_{B_x}, D_{B_x} )_{x\in B}$
is a GLS*$(1,m;2,(3, $ $K_0,\{3\},\{m+2\}),m|B|+2)$.
Then there is exactly one block $G_x\cup S\in \cup_{x\in X} {\cal F}_{x}$ from Step 3.

(2) $|T\cap (G_x\cup S)|\le 2$ for any $x\in X$. We distinguish 2 subcases.

(I) $|T\cap (G_x\cup S)|=2$ for some $x\in X$.
If $\infty\in \{y,z,g\}$, it is clear that only one element is $\infty$ and the other two are distinct.
Without loss of generality, let $g=\infty$ and $y=x$. Then $(g,c)\in S$ and $(y,a)\in G_x$.
There is exactly one block $B\in \A_1$ such that $\{x,z\}\subset B$ since $(X,\A_1)$ is an S$(2,K_1,u)$.
So $T$ is exactly contained in one block of  $\B_{B_x}$ since $(B\times\mathbb{Z}_m,\infty_{1}, \infty_{2}, \B_{B_x}, D_{B_x} )_{x\in B}$ is a GLS*.
Otherwise $\infty\not\in \{y,z,g\}$, then $|\{y,z,g\}|=2$.
Without loss of generality, let $g=z=x$ and $y\neq x$. Then  $(z,a),(g,c)\in G_x$.
There is exactly one block $B\in \A_1$ such that $\{x,y\}\subset B$ since $(X,\A_1)$ is an S$(2,K_1,u)$.
So $T$ is exactly contained in one block of  $\B_{B_x}$.

(II) $|T\cap (G_x\cup S)|\le 1$ for any $x\in X$.
Then then $|\{y,z,g\}|=3$ and $\infty\not\in \{y,z,g\}$.
There is exactly one block $B\in \A_1\cup {\cal T}$ such that $\{y,z,g\}\subset B$
since  $(X,\G,\A_1,{\cal T})$ is a 1-FG$(3,(K_1,K_{T}), u)$.
If $B\in \A_1$, by Step 1 $T$ is exactly contained in one block of $\cup_{x\in B}\B_{B_x}$ since
$(B\times\mathbb{Z}_m,\infty_{1}, \infty_{2}, \B_{B_x}, D_{B_x} )_{x\in B}$ is a GLS.
Otherwise $B\in {\cal T}$. There is a element $x\in B$ such that $x\not\in \{y,z,g\}$ since $|B|\ge 4$.
 By Step 2 $T$ is exactly contained in one block of $\C_B(x,i)$ since  $((B\backslash \{x\})\times \mathbb{Z}_m, \Gamma_{B}\backslash\{G_{x}\}, \C_B(x,i))$ is a $(3,1)$-GDD$(m^{|B|-1})$.

Thirdly, we show that each block $A$,  $A\in\F_{x}$ and $|A|\ge 4$, appears $m$ times in $\F_{x}$ and $|A|-2$ times in the multiset $\{\F_{x}:x\in X\}$, respectively.

Let $A\in\F_{x}$ and $|A|\ge 4$.
Then $A$ must come from Step 1.
Thus, there is a block  $B\in \A_1$ such that $A\in \cup_{x\in B}\B_{B_x}$.
 Since $(B\times \mathbb{Z}_m, \infty_{1}, \infty_{2}, {\cal B}_{B_x},D_{B_x})_{x\in X}$ is a
GLS*,  $A$ appears $m$ times in $\B_{B_x}$ and $|B|-2$ times in the multiset $\{\B_{B_x} : x\in X\}$, respectively.
So  $A$ appears $m$ times in $\F_x$ and $w(|B|-2)=|A|-2$ times in the multiset $\{\F_x:x\in X\}$ from Step 1, respectively.

Finally, we show that  each ordered pair $((x,a),(y,b))$ of distinct elements of $X\times \mathbb{Z}_m$,
not contained in some block $A$ with $\{\infty_{1},\infty_{2}\}\subset A$, appears in exactly one $D_{x}$.

For each ordered pair $(x, y)$ of distinct elements of $X\times \mathbb{Z}_m$,
not contained in some block $\{\infty_{1},\infty_{2}\}\in B$,
 there exists a unique block $B\in \A_1$ such that $\{x,y\} \subset A$ since $(X,\A_1)$ is an S$(2,K_1,u)$.
 Then by Step 1 $((x,a),(y,b))\in \cup_{x\in B} A(D_{B_x})$ since definition of GLS*.

Now, the proof is complete.\qed

\begin{theorem}\label{w=3^u} There exist a GLS*$(1,3;2,(3,\{3\},\{3\},\{5\}),3u+2)$ and a GLS$(1,3;2,(3,\{3\},$ $\{3\},\{3\}),3u+2)$
for any  $u\equiv 1, 3 \pmod 6$.
\end{theorem}

\noindent {\it Proof:}
There is an S$(3,4,u+1)$ in \cite{H1960}.
Delete a point from its point set to obtain a 1-FG$(3,(\{3\},\{4\}), u)$.
A GLS*$(1,3;2,(3, \{3\}, \{3\}, \{5\}), 11)$ exists by Lemma \ref{v=11}.
Apply Construction \ref{PCS} with $m=3$ to
get a GLS*$(1,3;2,(3,\{3\},\{3\},\{5\}),3u+2)$.
There is  a  GLS*$(1,3;2, (3, \{3\}, $ $\{3\}, \{3\} ), 5)$ by Lemma \ref{v=5}.
We apply Construction \ref{Dasui} to get a GLS$(1,3;2, (3, \{3\}, $ $\{3\}, \{3\} ), 3u +2)$.\qed

\section{ The spectrum for  $(3,\lambda)$-LGDD$(g^u)$}

In this section, we shall  prove our main result on $(3,\lambda)$-LGDD$(g^u)$.
We start with a direct construction.

\begin{lemma}\label{3^8}
There exists a $(3,2)$-LGDD$(3^8)$.
\end{lemma}

\noindent {\it Proof:}
 We shall construct a $(3,2)$-LGDD$(3^8)$ with point set $X=\mathbb{Z}_{24}$
and group set  $\G=\{ \{i,i+ 8,i+16\} : 0\le i\le7 \}$.
The block set $\B_i$ of the $i$-th  simple $(3,2)$-GDD$(3^8)$  can be generated from an initial block set
$\A_{i}$ by $+1 \pmod{24}$,   $0\le i\le 8$.

{\small
\begin{tabular}{lllllllllllllll}
 $\A_0 = $
&$\{( 0,   1,    2)$,
&$  (0,    2,    5)$,
&$  (0,    3,    9)$,
&$  (0,    4,   13)$,
&$  (0,    4,   14)$,
&$  (0,    5,   12)$,
&$  (0,    6,   13)$\},\\

$\A_1 = $
&$\{(0,    1,    3)$,
&$  (0,    1,    4)$,
&$  (0,    2,   11)$,
&$  (0,    4,   10)$,
&$  (0,    5,   14)$,
&$  (0,    5,   17)$,
&$  (0,    6,   17)$\},\\

$\A_2 = $
&$\{(0,    1,    5)$,
&$  (0,    1,    6)$,
&$  (0,    2,    6)$,
&$  (0,    2,   15)$,
&$  (0,    3,   12)$,
&$  (0,    3,   13)$,
&$  (0,    7,   14)$\},\\

$\A_3 = $
&$\{(0,    1,    7)$,
&$  (0,    1,   10)$,
&$  (0,    2,   12)$,
&$  (0,    2,   21)$,
&$  (0,    3,    7)$,
&$  (0,    4,   15)$,
&$  (0,    5,   11)$\},\\

$\A_4 = $
&$\{(0,    1,   11)$,
&$  (0,    1,   12)$,
&$  (0,    2,    7)$,
&$  (0,    2,    9)$,
&$  (0,    3,    6)$,
&$  (0,    4,    9)$,
&$  (0,    4,   18)$\},\\

$\A_5 = $
&$\{(0,    1,   13)$,
&$  (0,    1,   14)$,
&$  (0,    2,   17)$,
&$  (0,    2,   20)$,
&$  (0,    3,   18)$,
&$  (0,    3,   20)$,
&$  (0,    5,   10)$\},\\

$\A_6 = $
&$\{(0,    1,   15)$,
&$  (0,    1,   21)$,
&$  (0,    2,   13)$,
&$  (0,    2,   19)$,
&$  (0,    3,   10)$,
&$  (0,    4,   19)$,
&$  (0,    6,   12)$\},\\

$\A_7 = $
&$\{(0,    1,   18)$,
&$  (0,    1,   19)$,
&$  (0,    2,    4)$,
&$  (0,    3,   14)$,
&$  (0,    3,   15)$,
&$  (0,    4,   11)$,
&$  (0,    5,   15)$\},\\

$\A_8 = $
&$\{(0,    1,   20)$,
&$  (0,    1,   22)$,
&$  (0,    2,   14)$,
&$  (0,    3,   17)$,
&$  (0,    4,   17)$,
&$  (0,    5,   18)$,
&$  (0,    6,   15)$\}.
\end{tabular}}

Let $R=\{0,1,2,\ldots, 8\}$.
It is easy to check that $(X,\G, \B_r)_{r\in R}$ is a $(3,2)$-LGDD$(3^8)$.
\qed

\begin{construction}\label{LGDD-jia}{\rm(\cite{L1997})}
 If there is a $(3,\lambda)$-LGDD$(g^u)$,
 then there exists a $(3,\lambda)$-LGDD$(gm)^u$ for any $m\ge 1$.
\end{construction}

\begin{construction}\label{LGDD-hebing}
Let $t$ be a positive integer.
If there is a $(3,\lambda)$-LGDD$(g^u)$ and $g(u-2)\equiv 0 \pmod {\lambda t}$,
then there exists a $(3,\lambda t)$-LGDD$(g^u)$.
\end{construction}

\noindent {\it Proof:}
Suppose $\{(X,\G,\B_r): 1\le r\le \frac{g(u-2)}{\lambda}\}$ is a $(3,\lambda)$-LGDD$(g^u)$.
Let $${\cal A}_{i}=\bigcup\limits_{r=it+1}^{(i+1)t}\B_r, \ 0\le i\le \frac{g(u-2)}{\lambda t}-1.$$
Then $\{(X, \G,{\cal A}_i): 0\le i\le \frac{g(u-2)}{\lambda t}-1\}$ is a  $(3,\lambda t)$-LGDD$(g^u)$.\qed

\begin{lemma}\label{3^u}
There exists a $(3,2)$-LGDD$(3^u)$ for any $u\equiv 2 \pmod 6$ and $u\ge8$.
\end{lemma}

\noindent {\it Proof:} We distinguish 2 cases.

1. $u\equiv 2\pmod {12}$.
By Theorem \ref{(3,2)}
 we have an LS*$(1,2;2,(3,\{3,14,26,50,98\}),u)$.
An LS*$(1,2;2,(3,\{3,14\}),50)$
 and an LS*$(1,2;2,(3,\{3,26\}),98)$  exist by Lemma \ref{v=50} and Lemma \ref{v=98}, respectively.
Then we apply Construction \ref{LS-LGDD} (1) to get
 an LS*$(1,2;2,(3,\{3,14, 26\}),u)$.
Further, we apply Construction \ref{LS-LGDD} (2) with
 an LS$(1,2;2, $ $(3,\{3, 5\}), 14)$ and an LS$(1,2;2, $ $(3,\{3, 5\}), 26)$ from Lemma \ref{v=k}
to obtain
 an LS$(1,2;2,(3,\{3,5\}),u)$.
Finally, we apply Construction \ref{w=2} (3) to get a $(3,2)$-LGDD$(3^u)$.

2. $u\equiv 8\pmod {12}$. By  Theorem~\ref{(3,2)} an LS*$(1,2;2,(3,\{3,8\}),u)$ exists.
Apply Construction \ref{LS-LGDD} (4) with a $(3,2)$-LGDD$(3^8)$ from Lemma \ref{3^8} to get a $(3,2)$-LGDD$(3^u)$.
\qed

\begin{lemma}\label{2^u}
There exists a $(3,3)$-LGDD$(2^u)$ for any $u\equiv 2 \pmod 6$ and $u\ge8$.
\end{lemma}
\noindent {\it Proof:}
Assume that $u=2^{n}\cdot v+2$, where $v\equiv 3 \pmod 6$ and $n\ge 1$.
We distinguish 2 cases.

1.  $v\equiv 3,9 \pmod {18}$.
By Theorem \ref{w=3^u} there is an LS$(1,3;2, (3, \{3\},\{3\}, \{3\} ), v +2)$.
Apply Construction \ref{w=2} (2) to get an LS$(1,3; 2, (3, \{3,4\},\{3\}, \{2+2\}),2v+2)$.
Then we continue to use Construction \ref{w=2} (2) to get an LS$(1,3; 2, (3, \{3,4\},\{3\}, \{2^2+2\}),2^2\cdot v+2)$.
 Thus, by induction on $n$, we can prove that an LS$(1,3; 2, (3, \{3,4\},\{3\}, \{2^n+2\}),2^n\cdot v+2)$ exists.
 At last, we apply Construction \ref{w=2} (3) to get a $(3,3)$-LGDD$(2^u)$, where the input designs
 a $(3,1)$-LGDD$(2^4)$ and a $(3,1)$-LGDD$(2^{2^n+2})$  exist by Theorem \ref{LGDD}.

2.  $v\equiv 15 \pmod {18}$. Let $v=3w$. Then we have $w\equiv 5\pmod 6$.
Apply Construction \ref{w=w} with a GLS$(1,3;2, (3, \{3\},\{3\}, \{3\} ),5)$ from Lemma~\ref{v=5}
to get a GLS$(1,3;2, (3, \{3\},\{3\}, \{w+2\} ),3w +2)$.
Apply Construction \ref{w=2} (2) to get an LS$(1,3; 2, (3, \{3,4\},\{3\}, \{2w+2\}),2v+2)$.
Then we continue to use Construction \ref{w=2} (2) to get an LS$(1,3; 2, (3, \{3,4\},\{3\}, \{2^2\cdot w+2\}),2^2\cdot v+2)$.
 Thus, by induction on $n$, we can prove that an LS$(1,3; 2, (3, \{3,4\},\{3\}, \{2^n\cdot w+2\}),2^n\cdot v+2)$ exists.
 At last, we apply Construction \ref{w=2} (3) with a $(3,1)$-LGDD$(2^4)$ and a $(3,1)$-LGDD$(2^{2^n\cdot w+2})$  from Theorem \ref{LGDD}
 to get a $(3,3)$-LGDD$(2^u)$.
\qed

\begin{lemma}\label{lambda=1}
Let gcd$(\lambda,6)=1$ and $\lambda >1$. Then there exists a $(3,\lambda)$-LGDD$(g^u)$ if and only if $u\ge 3$, $2\le \lambda\le g(u-2)$,
$g(u-1)\equiv 0 \pmod 2$, $g^2u(u-1)\equiv 0 \pmod 6$ and $g(u-2)\equiv 0 \pmod {\lambda}$.
 \end{lemma}

\noindent {\it Proof:}
By Theorem \ref{Ne-LGDD} we only need to prove the sufficient condition.
For any $(g, u) \neq (1,7)$ there is a $(3,1)$-LGDD$(g^u)$ by Theorem~\ref{LGDD}.
Then we apply Construction \ref{LGDD-hebing} to get a $(3,\lambda)$-LGDD$(g^u)$.
For $(g, u) = (1,7)$,
since gcd$(\lambda,6)=1$, $\lambda >1$, and $\lambda\le g(u-2)$,
we have $\lambda=5$.
Then there is a $(3,5)$-LGDD$(1^7)$ which is also a simple $(3,5)$-GDD$(1^7)$ by Theorem \ref{SGDD}.
\qed

\begin{lemma}\label{lambda=2}
Let gcd$(\lambda,6)=2$.
Then there exists a $(3,\lambda)$-LGDD$(g^u)$
if and only if $u\ge 3$, $2\le \lambda\le g(u-2)$,
$g^2u(u-1)\equiv 0 \pmod 6$ and $g(u-2)\equiv 0 \pmod {\lambda}$.
\end{lemma}

\noindent {\it Proof:}
By Theorem \ref{Ne-LGDD} we only need to prove the sufficient condition.
Let $\lambda=2l$. Then we have $1\le l\le g(u-2)/2$.
We distinguish  3 cases as follows.

1.  gcd$(g,6)=1$. Then we have $u\equiv 0\pmod 2$ from $g(u-2)\equiv 0 \pmod 2$, and  $u\equiv 0,1 \pmod 3$ from $g^2u(u-1)\equiv 0 \pmod 6$.
So $u\equiv 0,4 \pmod 6$ and there is a $(3,2)$-LGDD$(1^u)$ by Theorem~\ref{LGDD}.
Apply Construction \ref{LGDD-jia}  with $\lambda=2$ and $m=g$ to get a $(3,2)$-LGDD$(g^u)$.
Further, we apply Construction \ref{LGDD-hebing} with $\lambda=2$ and $t=l$ to get a $(3,2l)$-LGDD$(g^u)$.

2.  gcd$(g,6)=2$ or $6$.
Then $g(u-1)\equiv 0 \pmod 2$ and $g^2u(u-1)\equiv 0 \pmod 6$.
So there exists a $(3,1)$-LGDD$(g^u)$ by Theorem~\ref{LGDD}.
Apply Construction \ref{LGDD-hebing} with $t=2l$ and $\lambda=1$ to
 obtain a $(3,2l)$-LGDD$(g^u)$.

3. gcd$(g,6)=3$.
Then $u\equiv 0 \pmod 2$.
If $u\equiv 0,4 \pmod 6$, we take a $(3,2)$-LGDD$(1^u)$ from  Theorem~\ref{LGDD} and use Construction \ref{LGDD-jia} with $\lambda=2$ and $m=g$ to get a $(3,2)$-LGDD$(g^u)$. If $u\equiv 2 \pmod 6$, we take a $(3,2)$-LGDD$(3^u)$ from  Lemma \ref{3^u} and use Construction \ref{LGDD-jia} with $\lambda=2$ and $m=g/3$ to get a $(3,2)$-LGDD$(g^u)$.
Further, we apply Construction \ref{LGDD-hebing} with $\lambda=2$ and $t=l$ to get a $(3,2l)$-LGDD$(g^u)$.
\qed

\begin{lemma}\label{lambda=3}
Let gcd$(\lambda,6)=3$.
Then there exists a $(3,\lambda)$-LGDD$(g^u)$ if and only if $u\ge 3$, $3\le \lambda\le g(u-2)$,
$g(u-1)\equiv 0 \pmod 2$ and $g(u-2)\equiv 0 \pmod {\lambda}$.
\end{lemma}

\noindent {\it Proof:}
By Theorem \ref{Ne-LGDD} we only need to prove the sufficient condition.
Let $\lambda=3l$. Then we have $1\le l\le g(u-2)/3$.
We distinguish  3 cases as follows.

1.  gcd$(g,6)=1$.  Then we have $u\equiv 1 \pmod 2$ from $g(u-1)\equiv 0 \pmod 2$, and $u\equiv 2 \pmod 3$ from $g(u-2)\equiv 0 \pmod {3}$.
So $u\equiv 5 \pmod 6$ and
there is a $(3,3)$-LGDD$(1^u)$ by Theorem~\ref{LGDD}.
Apply Construction \ref{LGDD-jia} with $\lambda=3$ and $m=g$ to get a $(3,3)$-LGDD$(g^u)$.
Further, we apply Construction \ref{LGDD-hebing}  with $\lambda=3$ and $t=l$ to obtain a $(3,3l)$-LGDD$(g^u)$.

2. gcd$(g,6)=2$.  Then $u\equiv 2 \pmod 3$. If $u\equiv 5 \pmod 6$, we take a $(3,3)$-LGDD$(1^u)$ from  Theorem~\ref{LGDD} and use Construction \ref{LGDD-jia} with $\lambda=3$ and $m=g$ to get a $(3,3)$-LGDD$(g^u)$. If $u\equiv 2 \pmod 6$, we take a $(3,3)$-LGDD$(2^u)$ from  Lemma \ref{2^u} and use Construction \ref{LGDD-jia} with $\lambda=3$ and $m=g/2$ to get a $(3,3)$-LGDD$(g^u)$.
Further, we apply Construction \ref{LGDD-hebing} with $\lambda=3$ and $t=l$ to get a $(3,3l)$-LGDD$(g^u)$.

3. gcd$(g,6)=3$ or $6$.  Then $g(u-1)\equiv 0 \pmod 2$ and $g^2u(u-1)\equiv 0 \pmod 6$.
So there exists a $(3,1)$-LGDD$(g^u)$ by Theorem~\ref{LGDD}.
Apply Construction \ref{LGDD-hebing} with $t=3l$ and $\lambda=1$ to get a $(3,3l)$-LGDD$(g^u)$.
\qed

\begin{lemma}\label{lambda=6}
Let gcd$(\lambda,6)=6$.
Then there exists a $(3,\lambda)$-LGDD$(g^u)$ if and only if
$u\ge 3$,   $6\le \lambda\le g(u-2)$ and $g(u-2)\equiv 0 \pmod {\lambda}$.
\end{lemma}

\noindent {\it Proof:}
By Theorem \ref{Ne-LGDD} we only need to prove the sufficient condition.
Let $\lambda=6l$. Then we have $1\le l\le g(u-2)/6$.
We distinguish  4 cases as follows.

1. gcd$(g,6)=1$. Then $u\equiv 2 \pmod 6$.
So there is a $(3,6)$-LGDD$(1^u)$ by Theorem~\ref{LGDD}.  Apply Construction \ref{LGDD-jia} with $\lambda=6$ and $m=g$ to get a $(3,6)$-LGDD$(g^u)$.
Further, we apply Construction \ref{LGDD-hebing} with $t=l$ and $\lambda=6$ to get a $(3,6l)$-LGDD$(g^u)$.

2.  gcd$(g,6)=2$. Then $g(u-1)\equiv 0 \pmod 2$ and $g(u-2)\equiv 0 \pmod {3}$.
So there is a $(3,3)$-LGDD$(g^u)$ by Lemma~\ref{lambda=3}.
Apply Construction \ref{LGDD-hebing} with $t=2l$ and $\lambda=3$ to get a $(3,6l)$-LGDD$(g^u)$.

3.  gcd$(g,6)=3$. Then $g(u-2)\equiv 0 \pmod {2}$ and $g^2u(u-1)\equiv 0 \pmod 6$.
So there is a $(3,2)$-LGDD$(g^u)$ by Lemma~\ref{lambda=2}.
Apply Construction \ref{LGDD-hebing} with $t=3l$ and $\lambda=2$ to get a $(3,6l)$-LGDD$(g^u)$.

4.  gcd$(g,6)=6$. Then  $g(u-1)\equiv 0 \pmod 2$ and $g^2u(u-1)\equiv 0 \pmod 6$.
So there is a $(3,1)$-LGDD$(g^u)$ by Theorem~\ref{LGDD}.
Apply Construction \ref{LGDD-hebing} with $t=6l$ and $\lambda=1$ to get a $(3,6l)$-LGDD$(g^u)$.
\qed

\noindent {\bf Proof of Theorem \ref{Main-L}:}
Combining Theorem~\ref{LGDD} and Lemmas \ref{lambda=1}-\ref{lambda=6},
we have come to the conclusion.
\qed

\section{ The spectrum for simple $(3,\lambda)$-GDD$(g^u)$}

In this section, we shall prove our main result for simple group divisible designs.

\begin{lemma}\label{la=1}
Let  gcd$(\lambda,6)=1$  and $g>1$. Then there exists a simple $(3,\lambda)$-GDD$(g^u)$ if and only if
$g(u-1)\equiv 0 \pmod 2$ and $g^2u(u-1)\equiv 0 \pmod 6$, $1\le \lambda\le g(u-2)$ and $u\ge 3$.
 \end{lemma}

\noindent {\it Proof:} Apply Lemma~\ref{LS-SGDD} with  a $(3,1)$-LGDD$(g^u)$ from Theorem~\ref{LGDD} to get a simple $(3,\lambda)$-GDD$(g^u)$.
\qed

\begin{lemma}\label{la=2}
Let  gcd$(\lambda,6)=2$  and $g>1$. Then there exists a simple $(3,\lambda)$-GDD$(g^u)$ if and only if
$g^2u(u-1)\equiv 0 \pmod 6$, $2\le \lambda\le g(u-2)$ and $u\ge 3$.
\end{lemma}

\noindent {\it Proof:}
Let $\lambda=2l$. Then we have $1\le l\le g(u-2)/2$.
By Theorem \ref{Ne-SGDD}, we only need to prove that the necessary condition is also sufficient. We distinguish 2 cases as follows.

1. $g(u-1)\equiv 0 \pmod 2$.
Apply Lemma~\ref{LS-SGDD} with $\lambda=1$, $t=2l$, and a $(3,1)$-LGDD$(g^u)$ from Theorem~\ref{LGDD} to get a simple $(3,2l)$-GDD$(g^u)$.

2.  $g(u-1)\equiv 1 \pmod 2$.
Then we have   $u\equiv 0 \pmod 2$.
So $g(u-2)\equiv 0 \pmod 2$.
By Lemma~\ref{lambda=2} there is a $(3,2)$-LGDD$(g^u)$.
Apply Lemma~\ref{LS-SGDD} with $\lambda=2$ and $t=l$ to get a simple $(3,2l)$-GDD$(g^u)$.
\qed

\begin{lemma}\label{la=3}
Let gcd$(\lambda,6)=3$ and $g>1$. Then there exists a simple $(3,\lambda)$-GDD$(g^u)$ if and only if $g(u-1)\equiv 0 \pmod 2$,
$3\le \lambda\le g(u-2)$ and $u\ge 3$.
\end{lemma}

\noindent {\it Proof:}
 Let $\lambda=3l$. Then we have $1\le l\le g(u-2)/3$.
By Theorem \ref{Ne-SGDD}, we only need to prove that the necessary condition is also sufficient. We distinguish 2 cases as follows.

1. $g^2u(u-1)\equiv 0 \pmod 6$. Apply Lemma~\ref{LS-SGDD} with $\lambda=1$, $t=3l$, and a $(3,1)$-LGDD$(g^u)$ from Theorem~\ref{LGDD} to get a simple $(3,3l)$-GDD$(g^u)$.

2. $g^2u(u-1)\equiv 2, 4\pmod 6$.
Then we have $u\equiv 2\pmod 3$.
So $g(u-2)\equiv 0 \pmod 3$.
By Lemma~\ref{lambda=3} there is a  $(3,3)$-LGDD$(g^u)$.
Apply Lemma~\ref{LS-SGDD} with $\lambda=3$ and $t=l$ to get a simple $(3,3l)$-GDD$(g^u)$.
\qed

\begin{lemma}\label{la=6}
Let gcd$(\lambda,6)=6$ and $g>1$. Then there exists a simple $(3,\lambda)$-GDD$(g^u)$ if and only if
$6\le \lambda\le g(u-2)$ and $u\ge 3$.
\end{lemma}

\noindent {\it Proof:}  Let $\lambda=6l$. Then we have $1\le l\le g(u-2)/6$.
By Theorem \ref{Ne-SGDD}, we only need to prove that the necessary condition is also sufficient. We distinguish 4 cases as follows.

1. $g(u-1)\equiv 0 \pmod 2$, $g^2u(u-1)\equiv 0 \pmod 6$.
Apply Lemma~\ref{LS-SGDD} with $\lambda=1$, $t=6l$, and a $(3,1)$-LGDD$(g^u)$ from Theorem~\ref{LGDD} to get a simple $(3,6l)$-GDD$(g^u)$.

2. $g(u-1)\equiv 0 \pmod 2$, $g^2u(u-1)\equiv 2, 4 \pmod 6$.
Then we have $u\equiv 2\pmod 3$.
So $g(u-2)\equiv 0 \pmod 3$.
By Lemma~\ref{lambda=3} there is a  $(3,3)$-LGDD$(g^u)$.
Apply Lemma~\ref{LS-SGDD} with $\lambda=3$ and $t=2l$ to get a simple $(3,6l)$-GDD$(g^u)$.

3. $g(u-1)\equiv 1 \pmod 2$, $g^2u(u-1)\equiv 0 \pmod 6$.
Then we have $u\equiv 0\pmod 2$.
So $g(u-2)\equiv 0 \pmod 2$.
By Lemma~\ref{lambda=2} there is a $(3,2)$-LGDD$(g^u)$.
Apply Lemma~\ref{LS-SGDD} with $\lambda=2$ and $t=3l$ to get a simple $(3,6l)$-GDD$(g^u)$.

4. $g(u-1)\equiv 1 \pmod 2$, $g^2u(u-1)\equiv 2, 4 \pmod 6$.
Then we have $u\equiv 2 \pmod 6$.  So $g(u-2)\equiv 0 \pmod 6$.
By Lemma~\ref{lambda=6} there is a $(3,6)$-LGDD$(g^u)$.
Apply Lemma~\ref{LS-SGDD} with $\lambda=6$ and $t=l$ to get a simple $(3,6l)$-GDD$(g^u)$.
\qed

\noindent {\bf Proof of Theorem \ref{Main-S}:}
Combining Theorem \ref{SGDD} and Lemmas~\ref{la=1}-\ref{la=6},
we have come to the conclusion.
\qed

\end{document}